\documentclass[11pt]{article}

\usepackage[margin=1in]{geometry}
\usepackage{amsmath,amsfonts,bm}
\usepackage{amssymb,amsthm}
\usepackage{booktabs}
\usepackage{bbding}
\usepackage{mathtools}
\usepackage{xurl}
\usepackage{float}
\usepackage{paralist}
\usepackage[misc]{ifsym}
\usepackage{graphicx}
\usepackage{epstopdf}
\usepackage[colorlinks=true,urlcolor=blue,citecolor=red,linkcolor=blue]{hyperref}
\allowdisplaybreaks

\newcommand{\R}{\mathbb{R}}
\newcommand{\bA}{\mathbf{A}}
\newcommand{\bB}{\mathbf{B}}
\newcommand{\bC}{\mathbf{C}}
\newcommand{\be}{\mathbf{e}}
\newcommand{\f}{\mathbf{f}}
\newcommand{\h}{\mathbf{h}}
\newcommand{\bK}{\mathbf{K}}
\newcommand{\bn}{\mathbf{n}}
\newcommand{\bQ}{\mathbf{Q}}
\newcommand{\bu}{\mathbf{u}}
\newcommand{\bx}{\mathbf{x}}
\newcommand{\by}{\mathbf{y}}
\newcommand{\bz}{\mathbf{z}}
\newcommand{\bLambda}{\bm\Lambda}

\newcommand{\btheta}{\bm\theta}

\newcommand{\bzeta}{\bm\zeta}
\newcommand{\mL}{\mathcal{L}}
\newcommand{\mP}{\mathcal{P}}
\newcommand{\mB}{\mathcal{B}}
\newcommand{\mH}{\mathcal{H}}
\newcommand{\abs}[1]{| #1 |}
\newcommand{\norm}[1]{| #1 |}
\newcommand{\inp}[1]{\langle #1 \rangle}
\DeclareMathOperator{\dif}{d\!}

\newcommand{\doititle}[1]{\textit{#1}}
\newcommand{\arXiv}[1]{\href{https://arxiv.org/abs/#1}{arXiv:#1}}

\newtheorem{theorem}{Theorem}[section]

\theoremstyle{definition}

\newtheorem{remark}[theorem]{Remark}

\title{Stability in Training PINNs for Stiff PDEs: Why Initial Conditions Matter}
\author{%
Baoli Hao$^1$, Chun Liu$^1$, Ulisses Braga-Neto$^2$, Lifan Wang$^3$, and Ming Zhong$^{4,*}$\\[0.6em]
\small $^1$Department of Applied Mathematics, Illinois Institute of Technology, Chicago, IL 60616, USA\\
\small $^2$Department of Electrical and Computing Engineering, Texas A\&M University, College Station, TX 77843, USA\\
\small $^3$Department of Physics and Astronomy, Texas A\&M University, College Station, TX 77843, USA\\
\small $^4$Department of Mathematics, University of Houston, Houston, TX 77204, USA\\[0.4em]
\small $^*$Corresponding author: \href{mailto:mzhong3@central.uh.edu}{mzhong3@central.uh.edu}
}
\date{}

\begin{document}
\maketitle

% \begin{center}
% \small\itshape
% Published in revised form in \emph{Foundations of Data Science} (2026),
% \href{https://doi.org/10.3934/fods.2026016}{doi:10.3934/fods.2026016}.
% This version is free to download for private research and study only.
% Redistribution, resale, and derivative use are not permitted.
% \end{center}

\begin{abstract}
Training physics-informed neural networks (PINNs) on stiff, time-dependent PDEs remains a fundamental challenge due to optimization instabilities and gradient pathologies. Through a series of rigorous ablation studies and Neural Tangent Kernel (NTK) analysis, we identify that the exact enforcement of initial conditions (ICs) is a decisive factor in stabilizing the training landscape. We present the first systematic ablation of two core strategies: hard initial-condition constrained transformation and self-adaptive loss weighting. Our findings demonstrate that embedding ICs directly into the network architecture provides an implicit time-marching effect, effectively reducing spectral bias and enabling the solution of highly stiff benchmarks, including sharp transitions and high-frequency coupled systems, primarily under periodic boundary conditions, with a Dirichlet extension reported as an additional robustness check. This work provides a scalable framework for developing reliable and physically-consistent neural solvers for complex mechanical systems.
\end{abstract}


\medskip

\vspace{0.5em}

\section{Introduction}

Time-dependent partial differential equations (PDEs) constitute the
mathematical foundation of numerous scientific and engineering disciplines,
serving as essential tools for modeling a diverse range of physical phenomena.
Their applications span from the intricate dynamics of fluid flow, as described
by the Navier--Stokes equations, to the complex interplay of chemical reactions
in biological systems, captured by reaction--diffusion
models~\cite{allen1979microscopic, bazant2017thermodynamic, burgers1948mathematical, cahn1958free, cassani2021belousov, de1991turing, fibich2015nonlinear, horstmann2013precipitation,hyman2014liquid,kato1987nonlinear,kevrekidis2001discrete,kudryashov1990exact,lee2014physical,michelson1986steady,miranville2017cahn,nishiura1999skeleton,
shen2010numerical,
takatori2015towards, tian2015electrochemically,whitham2011linear,zhabotinsky2007belousov}.
However, the solution of stiff time-dependent PDEs can be exceedingly
challenging to obtain and often requires substantial computational resources
and specialized numerical techniques, particularly due to multiscale effects
and strong nonlinearities that are central to many problems in computational
mechanics and engineering science.

Machine learning methodologies have demonstrated remarkable success across
various scientific and engineering fields, transforming areas such as protein
folding prediction~\cite{jumper2021highly}, drug discovery and
development~\cite{carracedo2021review}, and even the optimization of
fundamental algorithmic processes such as matrix-vector
multiplication~\cite{fawzi2022discovering}. However, traditional machine
learning approaches, which rely on large datasets, can be impractical in
scientific contexts where data acquisition is expensive, time-consuming, or
simply infeasible. Moreover, purely data-driven models often lack physical
fidelity and interpretability, as they may fail to capture the underlying
physical principles governing the systems they aim to describe.
Physics-informed machine learning (PIML), and in particular physics-informed
neural networks (PINNs), address these limitations by integrating physical
laws directly into the learning process. This significantly reduces the amount
of training data required. By embedding constraints derived from governing PDEs
into the neural network architecture, PINNs enable accurate predictions even
when data are sparse or
noisy~\cite{cai2021deepm, chen2022using, chen2021solving,jin2021nsfnets,
liu2019multi,rad2020theory,raissi2018hidden,raissi2019physics,wang2023pinns, yang2019physics, zhu2021machine}. Despite their promise, effectively
training PINNs for stiff time-dependent PDEs remains a major challenge. The
difficulty of propagating information from initial and boundary conditions into
the computational domain, combined with the need to balance multiple loss
terms, often results in convergence issues and suboptimal
solutions~\cite{coutinho2023physics, mcclenny2020self,wang2020understanding,wang2020and}.

From a numerical perspective, this challenge is intimately connected to the
ill-conditioning of the PDE solution operator with respect to initial and
boundary data. Unlike traditional numerical methods, which exactly satisfy
initial conditions at the start of the solution loop, PINNs enforce these
conditions only in a soft, least-squares sense, making the training landscape
highly sensitive to their degree of satisfaction. This ill-conditioning
manifests as spectral bias in the optimization, causing PINNs to struggle
with propagating information from initial conditions to subsequent time
steps, which is a phenomenon widely documented in the
literature~\cite{haitsiukevich2022improved, krishnapriyan2021characterizing, mcclenny2020self,
wang2022respecting,wight2020solving}.

This paper presents a complete and systematic study investigating the
importance of training components for PINNs in solving stiff time-dependent
PDEs. We do not claim the hard-constrained trial-function
construction itself as novel: enforcing initial and/or boundary conditions
exactly through an ansatz of the form $\tilde{u} = \psi + \phi\, u_{nn}$
dates back at least to~\cite{lagaris1998artificial}, and has since been
developed in various forms (see Section~\ref{sec:related}). Our contribution
lies instead in how this construction is used and analyzed for stiff
evolution problems. To the best of our knowledge, this is the first work
to conduct a rigorous ablation study specifically isolating the role of
exact initial-condition enforcement as the decisive
stabilizing mechanism for stiff time-dependent forward problems, and to
provide a theoretical justification for its necessity
via Neural Tangent Kernel (NTK) analysis. We focus specifically on
two training schemes: hard constraints and self-adaptive loss
weights~\cite{mcclenny2020self}. The hard-constraint scheme directly embeds
initial and boundary conditions into the neural network architecture, reducing
the multi-objective loss to a single-objective problem and mitigating spectral
bias~\cite{wang2020and}. Through our ablation study, we conclude that
the hard-constraint transformation is the most decisive factor in stabilizing
training, effectively mimicking the behavior of traditional time integrators
that begin from exact initial data. We further provide a rigorous NTK analysis
demonstrating how this transformation simplifies the training dynamics and
reduces spectral bias in a theoretically grounded manner. We emphasize that
hard constraints are not the only effective strategy for handling stiff PDEs,
and that this approach can be integrated with other advanced PINN techniques,
such as causal PINNs~\cite{wang2020and}, time-marching
PINNs~\cite{wight2020solving}, residual-based-attention (RBA)
PINNs~\cite{anagnostopoulos2024residual}, curriculum
training~\cite{krishnapriyan2021characterizing}, and co-training
PINNs~\cite{zhong2024label}, to achieve further synergistic improvements in accuracy and robustness.

The rest of the paper is structured as follows. We compare our methodology to
other related works in Section~\ref{sec:related}; in
Section~\ref{section:method}, we discuss the hard-constrained framework and
other training methodologies used in our ablation and comparison study; we
present a detailed discussion of the theoretical mechanism behind the
hard-constrained transformation in Section~\ref{sec:theory_full}; we present the
results  from five major categories of examples in
Section~\ref{section:example}; we conclude the paper in
Section~\ref{section:conclude}, pointing out several future directions for this
line of research.

\begin{figure}[h!]
\centering
\includegraphics[scale = 0.33]{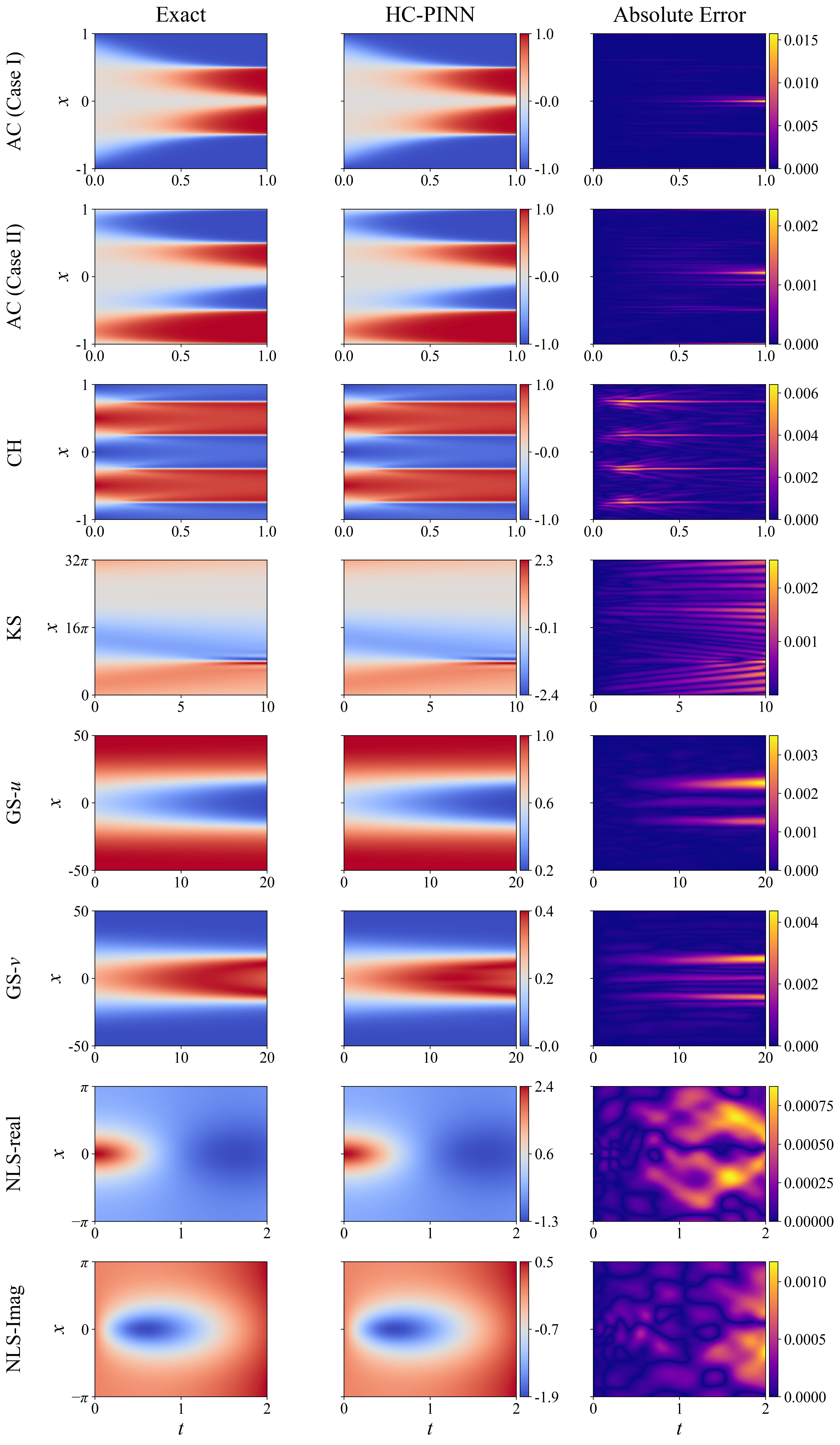}
\caption{Benchmark Results: Truth vs HC-PINN vs Absolute Error.}
\label{fig:all}
\end{figure}

\subsection{Related work}\label{sec:related}
PINNs have shown significant potential to solve a wide range of PDEs. However,
challenges arise, particularly when dealing with ``stiff'' PDEs, where explicit
numerical methods require extremely small time steps for stability. A primary
concern is the imbalance between the fitting of initial data and residual
components within the PINN loss
function~\cite{shin2020convergence,wang2020understanding,
wang2020and,wight2020solving}. This imbalance often leads gradient descent to prioritize
minimizing the residual loss over the data fitting loss, hindering convergence
to accurate solutions. This issue is exacerbated in forward problems, where
data is primarily concentrated at initial and boundary conditions, leaving
limited information within the PDE domain. Specifically, for time-dependent
problems, PINNs struggle to propagate information from initial conditions to
subsequent time steps, a phenomenon widely documented in the
literature~\cite{haitsiukevich2022improved,krishnapriyan2021characterizing,mcclenny2020self,
wang2022respecting,wight2020solving}.
Several strategies have been proposed to address these challenges. One
approach, introduced by~\cite{wight2020solving}, involves time-marching, where
the time domain is divided into smaller intervals, with PINNs being
sequentially trained on each interval, starting from the initial condition.
While effective, this method can be computationally expensive due to the need
to train multiple networks. To improve information propagation, researchers
have emphasized the importance of weighting initial condition data and residual
points near the initial time~\cite{mcclenny2020self,wight2020solving}.
\cite{mcclenny2020self} demonstrated that their self-adaptive PINNs can
autonomously learn these weights during training, particularly in scenarios
with complex initial conditions, such as the wave equation. Further
advancements include the ``sequence-to-sequence''
approach~\cite{krishnapriyan2021characterizing}, causality
PINN~\cite{wang2022respecting}, artificial viscosity
PINN~\cite{coutinho2023physics}, and co-training
PINN~\cite{zhong2024label}.
Beyond loss weighting and time marching, structural modifications to PINNs
have also been explored. \cite{braga2022characteristics} introduced
Characteristic-Informed Neural Networks (CINNs) for transport equations, where
PDE information is integrated into the neural network architecture. A similar
approach, albeit for Dirichlet boundary conditions, was proposed
in~\cite{lagaris1998artificial}. More recently, a growing
body of work has developed hard-constraint PINN formulations along these
architectural lines. \cite{lu2021physics} proposed hPINNs, which enforce
constraints exactly through trial-function transformations and the penalty/
augmented-Lagrangian method, primarily in the context of inverse design and
PDE-constrained optimization. \cite{sukumar2022exact} constructed approximate
distance functions to enforce Dirichlet, Neumann, and Robin boundary
conditions exactly, and related geometry-aware constructions have been used
to impose exact boundary conditions in PINN and variational (VPINN) settings.
Unified hard-constraint frameworks have also been proposed to handle multiple
boundary-condition types within a single formulation~\cite{liu2022unified}.
These works demonstrate the value of exact constraint enforcement, but they are
predominantly boundary-condition focused or aimed at inverse-design problems.
They do not isolate exact initial-condition enforcement as the stabilizing
mechanism for stiff time-dependent forward problems, nor provide an
ablation- and NTK-based explanation of this effect. Our work addresses this gap
by focusing specifically on exact initial-condition enforcement in stiff
time-dependent forward problems.

\section{Methodology}\label{section:method}
We consider the following setup for general time-dependent PDEs.  Let $u$ be an unknown function defined on $[0, T]\times\bar\Omega$ with $\Omega \subset \R^d$.  Here the physical domain $\Omega$ comes with a Lipschitz boundary $\partial\Omega$ and $\bar{\Omega} = \Omega \cup \partial\Omega$.  And we denote $\bx = \begin{pmatrix} x_1 & \cdots & x_d \end{pmatrix}^\top \in \Omega$. Then we say that $u$ is a solution of a time-dependent PDE, if $u$ satisfies the following
\begin{equation}\label{eq:pde_setup}
\begin{cases}
u_t(t, \bx) + \mP[u](t, \bx) &= f(t, \bx), \quad (t, \bx) \in (0, T]\times\Omega, \\
u(0, \bx) &= u_0(\bx), \quad \bx \in \bar\Omega, \\
\mB[u](t, \bx) &= g(t, \bx), \quad (t, \bx) \in [0, T]\times\partial\Omega, \\
\end{cases}.
\end{equation}
Here the operator $\mP$ has order $K$ and it is defined as
\[
\mP[u] = h(u, \{\partial^k_{x_1^{j_1}x_2^{j_2}\cdots x_d^{j_d}}u\}_{k = 1}^K), \quad \text{$h$ is any multivariate function},
\]
where $\partial^k_{x_1^{j_1}x_2^{j_2}\cdots x_d^{j_d}}u$ is a set that contains all of partial derivatives of $u$ with respective  to $\bx$ of $k^{th}$ order in the multivariate sense, i.e., for the powers of each coordinate, we have
\[
j_1, j_2, \cdots, j_d \ge 0 \quad \text{and} \quad j_1 + j_2 + \cdots + j_d = k.
\]
Furthermore, $\mathcal{B}$ denotes a boundary operator defined on
$[0,T]\times\partial\Omega$. We assume that the operator $\mathcal{P}$,
the boundary operator $\mathcal{B}$, and the prescribed data $u_0$, $f$,
and $g$ satisfy suitable conditions such that the corresponding initial-boundary
value problem is well posed. We also assume that the initial and boundary data
are compatible at $t=0$, namely,
\[
    \mathcal{B}[u_0](\bx)=g(0,\bx),
    \qquad \bx\in\partial\Omega.
\]
For Dirichlet boundary conditions, this reduces to
$u_0(\bx)=g(0,\bx)$ on $\partial\Omega$.
% Furthermore $\mathcal{B}$ is an operator defined for $(t, \bx)$ on the boundary.  We also assume that the functions $u_0$, $f$, and $g$ are user inputs that satisfy the desired regularity so that the existence and uniqueness of solutions for such a PDE are guaranteed.  We also assume that the compatibility condition, $g(0, \bx) = u_0(\bx)$ for $\bx \in \partial\Omega$, is satisfied to prevent any ill-conditioning of the solution.

\textbf{Types of Boundary Conditions}: one can consider many kinds of boundary conditions as follows
\[
\mathcal{B}[u](t, \bx;a, b) = au(t, \bx) + b\frac{\partial u}{\partial \bn}(t, \bx), \quad \text{Robin type},
\]
where $\bn$ is the outward normal vector to $\partial\Omega$ and $a, b$ are two fixed constants.  Special values of $a$ and $b$ can give two other BC cases, such as when $b = 0$ gives $\mathcal{B}[u](t, \bx) = u(t, \bx)$ (Dirichlet) and $a = 0$ gives the $\mathcal{B}[u](t, \bx) = \frac{\partial u}{\partial \bn}(t, \bx)$ (Neumann). The periodic boundary condition gives $d$ different equations, and they are,
\begin{equation}
\mathcal{B}[u](t, \bx) = u(t, \bx) - u(t, \bx + P\be_i) = 0, \quad i = 1, \cdots, d.
\end{equation}
Here $P$ is the period and $\be_i$ is the $i^{th}$ standard basis vector in $\R^d$, that is,
\[
\be_i =
\begin{bmatrix}
0 & \cdots & 0 & 1 & 0 & \cdots & 0
\end{bmatrix}^{\top},
\]
where the entry \(1\) is in the \(i\)-th position.
One can even consider a mixture of two or all of the above BC conditions.  For example, we may decompose the boundary as
\[
    \partial\Omega=\Gamma_1\cup\cdots\cup\Gamma_K,
\]
where each $\Gamma_i$ is associated with a different type of boundary condition.
% For example, we can have $\Omega = \Gamma_1 \cup \cdots \cup \Gamma_K$ with each $\Gamma_i$ defines a different kind of BC.
In this paper, we will mainly focus on the periodic type of boundary conditions.

Recent advancements in scientific machine learning integrate the principles of Physics, particularly the Physics-derived PDEs, into the training process of machine learning models. This integration enables the development of PINNs to solve for $u$ as follows: find an approximate solution from a set of deep neural networks $\mH_{\text{NN}}$, where each network has the same depth, the same number of neurons on each hidden layer, and the same activation functions on each layer, that minimizes the following loss functional
    \begin{equation}
    \text{Loss}(u_{nn}) = \text{Data Loss}(u_{nn}) + \lambda \cdot \underbrace{\text{IC Loss}(u_{nn}) + \text{BC Loss}(u_{nn}) + \text{PDE Loss}(u_{nn})}_{\text{Physical Loss}} \label{eq:loss-s}
\end{equation}
where $\lambda$ is a regularization parameter, and the loss components are defined as follows
    \[
\begin{cases}
    &\text{Data Loss}(u_{nn}) = \frac{1}{N_{Data}}\sum_{i = 1}^{N_{Data}}\abs{u_{nn}(t^{Data}, \bx_i^{Data}) - u(t^{Data}, \bx_i^{Data})}^2\\[0.2cm]
    &\text{IC Loss}(u_{nn}) = \frac{1}{N_{IC}}\sum_{i = 1}^{N_{IC}}\abs{u_{nn}(0, \bx_i^{IC}) - u_0(\bx_i^{IC})}^2,\\[0.2cm]
    &\text{BC Loss}_{\text{per}}(u_{nn})\! = \sum_{i = 1}^d\frac{1}{N_{BC}}\sum_{j \!= 1}^{N_{BC}}\abs{u_{nn}(t_j^{BC}, \bx_j^{BC}) \!- u_{nn}(t_j^{BC}, \bx_j^{BC} + P\be_i)}^2, \\
    & \quad \text{or}\\[0.2cm]
    &\text{BC Loss}_{\text{gen}}(u_{nn}) = \frac{1}{N_{BC}}\sum_{j = 1}^{N_{BC}}\abs{(\mB[u_{nn}] - g)(t_j^{BC}, \bx_j^{BC})}^2,\\[0.2cm]
    &\text{PDE Loss}(u_{nn}) = \frac{1}{N_{CL}}\sum_{i = 1}^{N_{CL}}\abs{(\frac{\partial u_{\rm nn}}{\partial t} + \mP[u_{nn}] - f)(t_i^{CL}, \bx_i^{CL})}^2
\end{cases}
\]
for $ u_{nn} \in \mathcal{H}_{\text{NN}}$. Note that $\text{BC Loss}_{\text{per}}(u_{nn})$ enforces continuity and periodicity along each spatial direction~$i$, while in $\text{BC Loss}_{\text{gen}}(u_{nn})$, $\mathcal{B}$ and $g$ denote the boundary operator (e.g., identity or normal derivative) and the prescribed boundary condition, respectively.  Here  $\{(t^{CL}, \bx^{CL})_i\}_{i = 1}^{N_{CL}} \in (0, T]\times\Omega$ are called collocation points, $\{(0, \bx^{IC})_i\}_{i = 1}^{N_{IC}} \in \{0\}\times\bar\Omega$ are the initial condition points, $\{(t^{BC}, \bx^{BC})_i\}_{i = 1}^{N_{BC}} \in [0, T]\times\partial\Omega$ are the boundary condition points, and $\{(t^{Data}, \bx^{Data})_i\}_{i = 1}^{N_{Data}} \in (0, T]\times\Omega$ are known data points (they can be noisy).  The minimizer, denoted as $u_{NN}$, will be an approximate solution to \eqref{eq:pde_setup}.

\subsection{Hard constraints}
The main motivation for considering a hard-constraint transformation for PINNs is the possible ill-conditioning of the PDE solution operator's dependence on IC and BC.  Traditional numerical methods do not encounter this kind of difficulty, as the IC data is exactly satisfied at the starting point of the solution loop.  PINNs, on the other hand, use $L^2$-loss to fit the IC/BC and PDE residual data, which causes the solutions to be highly sensitive to the tightness of fit for the IC in the training, as pointed out in~\cite{mcclenny2020self}.  By transforming the need of fitting IC/BC conditions out of PINNs, it reduces the spectral bias in the training.  Therefore, we consider the following transformation
\begin{equation}
    \tilde{u}(t, \bx) = \psi(t, \bx) + \phi(t, \bx)u_{nn}(t, \bx), \label{eq:trans}
\end{equation}
% Here the functions $\psi$ and $\phi$ are smooth functions with the following property
% \begin{equation}\label{eq:HC_condition}
% \begin{cases}
%  &\psi(0, \bx) = u_0(\bx), \quad \phi(0, \bx) = 0, \quad \text{for $\bx \in \bar\Omega$},\\
%  &\psi(t, \bx) = g(t, \bx), \quad \phi(t, \bx) = 0, \quad \text{for $(t, \bx) \in [0, T]\times\partial\Omega$}.
% \end{cases}
% \end{equation}
% In the case of periodic BC, we simply require $\psi$ and $\phi$ to be both $P$-periodic.
Here $\psi$ and $\phi$ are smooth functions chosen to enforce the initial
condition exactly. In particular, we require
\begin{equation}\label{eq:HC_condition}
    \psi(0,\bx)=u_0(\bx),
    \qquad
    \phi(0,\bx)=0,
    \qquad
    \bx\in\bar\Omega.
\end{equation}
For Dirichlet boundary conditions, the boundary condition can also be enforced
exactly by requiring
\begin{equation}\label{eq:HC_Dirichlet_condition}
    \psi(t,\bx)=g(t,\bx),
    \qquad
    \phi(t,\bx)=0,
    \qquad
    (t,\bx)\in[0,T]\times\partial\Omega.
\end{equation}
For periodic boundary conditions, instead of
\eqref{eq:HC_Dirichlet_condition}, we assume that $u_0$, $\psi$, and $\phi$
are $P$-periodic in each spatial direction. Together with the periodic neural
network construction introduced below, this guarantees that
$\widetilde{u}=\psi+\phi u_{\rm nn}$ satisfies the periodic boundary
conditions exactly.
When such transformation is used, then the training for $u^{nn}$ is updated as
\begin{equation}
    \text{Loss}_{HC}(u_{nn}) = \frac{1}{N_{CL}}\sum_{i = 1}^{N_{CL}}\abs{(\tilde{u}_t + \mP[\tilde{u}] - f)(t_i^{CL}, \bx_i^{CL})}^2, \quad \tilde{u} = \psi +\phi u_{nn}. \label{eq:loss-h}
\end{equation}
We have successfully reduced the multi-objective loss to a single-objective; however the information about IC and BC has been moved into the PDE residual loss through $\psi$ and $\phi$.

\textbf{Periodic Neural Networks}: even though we require $\psi$ and $\phi$ to be periodic, we still need to build the information of periodicity into $u_{nn}$, hence we consider the following composition
\[
u_{nn}(t, x) = f_{nn}(v(t, x)), \quad \text{where $f_{nn}$ is a neural network and $x \in \R$}.
\]
Here $v(t, x) = \begin{bmatrix} t  &  \mathcal{F}_m(x)\end{bmatrix}^\top$ with
\[
\mathcal{F}_m(x) = \begin{bmatrix}1 & \cos\left(\frac{2\pi}{P}x\right) &\sin\left(\frac{2\pi}{P}x\right) & \cdots  &\cos\left(\frac{2\pi}{P}mx\right) &  \sin\left(\frac{2\pi}{P}mx\right)\end{bmatrix}
\]
where $m$ is a positive integer hyperparameter controlling the number of Fourier modes. With this construction, each trigonometric feature is exactly $P$–periodic, and therefore any linear combination (and hence $u_{nn}$) will inherit $P$–periodicity~\cite{dong2021method}.  For $\bx = \begin{bmatrix}x & y\end{bmatrix}$, we can use the following
\[
v(t, x, y) = \begin{bmatrix} t & \mathcal{F}_m(x) & \mathcal{F}_m(y) \end{bmatrix};
\]
similarly for higher dimensional $\bx = \begin{bmatrix} x_1 & \cdots & x_d\end{bmatrix}$, we have
\[
v(t, \bx) = \begin{bmatrix} t & \mathcal{F}_m(x_1) & \cdots & \mathcal{F}_m(x_d) \end{bmatrix}.
\]

%
% \subsection{Other Training Enhancements}
% %
% We also consider mini-batching (to handle large sample size) and self-adaptive loss weights~\cite{mcclenny2020self}, please see sections \ref{sec:minibatch} and \ref{sec:SAPINN} for details.

\subsection{Mini-batching}\label{sec:minibatch}
We have eliminated the IC/BC losses, but the transformed PDE loss might create additional difficulties in training due to the introduction of IC/BC functions into the PDE.  We adopt a machine learning technique, called mini-batching, commonly used to reduce training complexity by evaluating the loss on randomly selected subsets of collocation data. In PINNs, this reduces memory costs, especially when handling large numbers of collocation points. Sampling mini-batches from collocation points also helps stabilize and accelerate optimization. So we consider mini-batching in the training for the PDE residual loss.  As pointed out in~\cite{wight2020solving}, the mini-batching technique is similar to time-marching resampling but without an explicit time grid.
\subsection{Self-adaptive PINNs}\label{sec:SAPINN}
The proper choice of penalty term $\lambda$ can be updated from derived formula as mentioned in~\cite{wang2020and}.  However, data-driven discovery of $\lambda$ is possible.  Following~\cite{mcclenny2020self}, the loss is further modified as
\[
\begin{aligned}
    &\text{IC/BC Loss}(\{\lambda^{IC}_i\}_{i = 1}^{N_{IC}}, \{\lambda^{BC}_i\}_{i = 1}^{N_{BC}})\\& = \frac{1}{N_{IC}}\sum_{i = 1}^{N_{IC}}F_{mask}(\lambda^{IC}_i)\abs{u_{nn}(0, \bx_i^{IC}) - u_0(\bx_i^{IC})}^2 \\
    &\quad+ \sum_{i = 1}^d\frac{1}{N_{BC}}\sum_{j = 1}^{N_{BC}}F_{mask}(\lambda^{BC}_j)\abs{u_{nn}(t_j^{BC}, \bx_j^{BC})  - u_{nn}(t_j^{BC}, \bx_j^{BC} + P\be_i)}^2,
\end{aligned}
\]
and
   \[
\text{Physical Loss}(\{\lambda^{CL}_i\}_{i = 1}^{N_{CL}})= \! \frac{1}{N_{CL}}\sum_{i = 1}^{N_{CL}}F_{mask}(\lambda^{CL}_i)\abs{(\frac{\partial u_{nn}}{\partial t}\! +\mathcal{P}(u_{nn}) - f)(t_i^{CL}, \bx_i^{CL})}^2.
\]
The set of loss weights, $\{\lambda^{IC}_i\}_{i = 1}^{N_{IC}}\cup\{\lambda^{BC}_i\}_{i = 1}^{N_{BC}}\cup\{\lambda^{CL}_i\}_{i = 1}^{N_{CL}}$, are designed in a way to increase when needed during any training epoch in order to combat the spectral bias.  Adaptive methods are essential to ensure that the neural network accurately addresses the challenging spots in the solutions of ``stiff'' PDEs. In our experiments, we employ self-adaptive PINNs as proposed by \cite{mcclenny2020self} to train PINNs adaptively. This approach involves fully trainable adaptation weights applied individually to each training point, allowing the neural network to independently identify and focus on the challenging regions of the solution.
\subsection{Performance measures}
Let $\{x_i,t_i\}_{i=1}^N$ be a collection of $N$ data points where we evaluate both the reference solution $u(x,t)$ and the neural network output $\mathcal{U}(x,t)$. To measure the accuracy of the trained model, we compute the relative $L^1$ norm $\mathcal{E}_1$ and relative $L^2$ norm $\mathcal{E}_2$,  defined by:
\begin{align*}
    \mathcal{E}_1 = \frac{\sum_{i=1}^N |\mathcal{U}(x_i,t_i)-u(x_i,t_i)|}{\sum_{i=1}^N |u(x_i,t_i)|},
\quad  \mathcal{E}_2 = \frac{\sqrt{\sum_{i=1}^N |\mathcal{U}(x_i,t_i)-u(x_i,t_i)|^2}}{\sqrt{\sum_{i=1}^N |u(x_i,t_i)|^2}}.
\end{align*}
These metrics quantify the discrepancy between the predicted and reference solutions in
the relative absolute-error ($L^1$) and quadratic/energy-error ($L^2$) senses, respectively.
The relative normalization enables fair comparisons across different scales and variables. To rigorously assess the accuracy of our results, we implemented a semi-implicit spectral method in Matlab. This allows us to make precise comparisons between the ground-truth solutions $u(x_i,t_i)$ and the PINN-generated solutions $\mathcal{U}(x_i,t_i)$.
\section{Theoretical foundation for hard constraints}\label{sec:theory_full}
This section provides the theoretical justification for the stability and conditioning advantages of the hard-constraint transformation in PINNs. While the transformation enforces both boundary and initial conditions, our emphasis is on the initial condition, as the focus of this paper is on solving stiff time-dependent PDEs with PINNs. Stiff time-dependent PDEs exhibit features across widely varying time scales, making them difficult to capture even for explicit time-integration methods in traditional numerical approaches.

By incorporating the initial condition directly into the trial function, the
hard-constraint transformation removes the initial-condition mismatch from the
optimization problem. The network is therefore trained only through the
transformed PDE residual, while every admissible approximation starts from the
prescribed initial state. This is analogous to conventional time-dependent
solvers, which evolve the solution from given initial data, and motivates the
improved training stability observed below.

% By removing the initial condition through the hard-constraint transformation, we obtain the following advantages: (i) we implicitly enforce the time direction during PINN training, reducing spectral bias arising from multi-objective loss formulations; (ii) we mimic the behavior of traditional numerical methods, which begin time integration from exact or approximated initial conditions; and (iii) the hard-constraint transformation behaves similarly to an implicit time integrator.

We divide the analysis into three parts: $(I)$  consistency and non-degeneracy of the hard-constraint transformation, $(II)$  Hessian and NTK characterization, and $(III)$  implications for training dynamics, including decay and spectral bias.

\subsection{Consistency and non-degeneracy of the hard-constraint transformation}
\label{sec:theory}

We start with the minimal regularity conditions on the pair $(\psi,\phi)$
under which the hard-constraint transformation is consistent with the original
initial-value problem and non-degenerate at the initial time.

\begin{theorem}[Consistency and Non-Degeneracy of the Transformation]
\label{thm:wellposed}
Assume that
\[
    \phi,\psi \in C^1\bigl([0,T];C^k(\bar{\Omega})\bigr),
\]
and that
\[
    \psi(0,\bx)=u_0(\bx),
    \qquad
    \phi(0,\bx)=0,
    \qquad
    \phi_t(0,\bx)\neq 0,
    \qquad
    \bx\in\bar{\Omega}.
\]
Let
\[
    u(t,\bx)=\psi(t,\bx)+\phi(t,\bx)u_{nn}(t,\bx).
\]
Then $u$ automatically satisfies the prescribed initial condition:
\[
    u(0,\bx)=u_0(\bx),
    \qquad
    \bx\in\bar{\Omega}.
\]

Moreover, $u$ satisfies the original PDE
\[
    \partial_tu+\mP[u]=f,
    \qquad
    (t,\bx)\in(0,T]\times\Omega,
\]
if and only if $u_{nn}$ satisfies the transformed PDE
\[
    \phi\partial_tu_{nn}
    +
    \mP[\psi+\phi u_{nn}]
    +
    \phi_tu_{nn}
    =
    f-\psi_t,
    \qquad
    (t,\bx)\in(0,T]\times\Omega.
\]
At the initial time, the value of the transformed variable is determined by
\[
    u_{nn}(0,\bx)
    =
    \frac{
        u_t(0,\bx)-\psi_t(0,\bx)
    }{
        \phi_t(0,\bx)
    },
    \qquad
    \bx\in\bar{\Omega}.
\]
Therefore, the condition $\phi_t(0,\bx)\neq 0$ prevents degeneracy in the
recovery of $u_{nn}(0,\bx)$ from the transformed representation.
\end{theorem}

\begin{remark}
Theorem~\ref{thm:wellposed} shows that the hard-constraint transformation
places the neural approximation on the correct initial-data manifold and
provides a non-degenerate representation at $t=0$. Intuitively, $\phi$
scales the learned correction while $\psi$ encodes the prescribed initial
condition. A bounded factor such as $1-e^{-Ct}$ or $t/T$ avoids the
large-$T$ growth associated with an unbounded choice such as $\phi=t$.

For example, one may choose
\[
    \psi(t,\bx)=e^{-Ct}u_0(\bx),
    \qquad
    \phi(t,\bx)=1-e^{-Ct}.
\]
In our implementation, we adopt the closely related choice
\[
    \psi(t,\bx)=e^{-Ct}u_0(\bx),
    \qquad
    \phi(t)=\frac{t}{T},
\]
which is bounded on $[0,T]$ for any $T$ and satisfies
\[
    \phi_t(0,\bx)=\frac{1}{T}\neq 0.
\]
This choice provides a simple normalized-time scaling while enforcing the
initial condition exactly.

When $\mP$ is linear and $\phi(t,\bx)=\phi(t)$, the transformed PDE takes
the simpler form
\[
    \phi\partial_t u_{nn}
    +
    \phi\mP[u_{nn}]
    +
    \phi' u_{nn}
    =
    f-\psi_t-\mP[\psi].
\]
If $\phi(t)\neq 0$ for $t>0$, then for positive times this equation can be
written as
\[
    \partial_tu_{nn}
    +
    \mP[u_{nn}]
    +
    \frac{\phi'}{\phi}u_{nn}
    =
    \frac{
        f-\psi_t-\mP[\psi]
    }{
        \phi
    }.
\]
Thus, $u_{nn}$ satisfies an equation of the same general form as the original
PDE, with an additional reaction term
$\frac{\phi'}{\phi}u_{nn}$ and a modified forcing term that carries the
information encoded in $\psi$, in particular the prescribed initial data.
\end{remark}

\subsection{NTK Perspective and training decay}
Finally, we provide a Neural Tangent Kernel (NTK) perspective on how hard constraints would work to reduce the spectral bias and training difficulties of PINNs on solving stiff time-dependent PDEs.  To simplify the discussion, we will focus on the PDEs with periodic boundary. Let \(\{\bz_i = (t_i, \bx_i)\}_{i=1}^N\subset(0,T]\times\Omega\) denote collocation (residual) points and \(\{\bzeta_j=(0, \bx_j)\}_{j=1}^M\) denote initial-condition sample points.  We define the two residual functions $r_s$ and $r_h$ as
\[
r_s(\bz_i) = \partial_t u_s(\bz_i) + \mP[u_s](\bz_i) - f(\bz_i),
\]
and
\[
r_h(\bz_i) = \phi(\bz_i)\partial_t u_h(\bz_i)
+ \mP[\psi + \phi u_h](\bz_i)
+ \phi_t(\bz_i)u_h(\bz_i)
- f(\bz_i) + \psi_t(\bz_i),
\]
where $u_s = u_s(t, \bx; \btheta_s)$ (soft constrained) and $u_h = u_h(t, \bx; \btheta_h)$ (hard constrained) are the two PINN approximations, which are automatically periodic due to the Fourier layer. The corresponding losses to train $u_s$ and $u_h$ are
\[
\text{Loss}_s(u_s(\btheta_s)) = \text{Loss}_{IC}(u_s(\btheta_s)) + \text{Loss}_{PDE}(u_s(\btheta_s))
\]
with
\[
\begin{cases}
\text{Loss}_{IC}(u_s(\btheta_s)) &= \frac{1}{2M}\sum_{i = 1}^{M}|u_s(\bzeta_i; \btheta_s) - u_0(\bx^{IC}_i)|^2, \quad \bzeta_i = (0, \bx^{IC}_i), \\[0.2cm]
\text{Loss}_{PDE}(u_s(\btheta_s)) &= \frac{1}{2N}\sum_{i = 1}^{N}|r_s(\bz_i; \btheta_s)|^2,
\end{cases}
\]
and  $\quad \quad
\text{Loss}_h(u_h(\btheta_h)) = \frac{1}{2N}\sum_{i = 1}^{N}|r_h(\bz_i; \btheta_h)|^2.
$

Both Losses can induce a continuous-time (in terms of $\tau$ not $t$ as $t$ represents the physical time) gradient flow systems as
\[
\frac{\dif\btheta_s(\tau)}{\dif\tau} = -\nabla\text{Loss}_s(u_s(\btheta_s(\tau))) \quad \text{and} \quad \frac{\dif\btheta_h(\tau)}{\dif\tau} = -\nabla\text{Loss}_h(u_h(\btheta_h(\tau))),
\]
with $\btheta_s(0) = \btheta_h(0) = \btheta_0$. Let the operators $\mL_s[u] \coloneqq \partial_tu + \mP[u]$ and $\mL_h[u] \coloneqq \phi\partial_tu + \mP[\psi + \phi u] + \phi_tu$.  For notational simplicity, in the following NTK dynamics we suppress the
normalization factors arising from $1/M$ and $1/N$, since they only rescale
the corresponding kernel contributions and do not affect the comparison
between the coupled soft-constrained system and the residual-only
hard-constrained system. Then with NTK, we obtain
\[
\begin{aligned}
\begin{bmatrix} \frac{\dif u_s(\bzeta; \btheta_s(\tau))}{\dif\tau} \\[0.2cm] \frac{\dif \mL_s[u_s](\bz; \btheta_s(\tau))}{\dif\tau}\end{bmatrix} = -\begin{bmatrix} \bK^s_{11}(\tau) & \bK^s_{12}(\tau) \\[0.2cm] \bK^s_{21}(\tau) & \bK^s_{22}(\tau)\end{bmatrix}\cdot\begin{bmatrix}u_s(\bzeta; \btheta_s(\tau)) - u_0(\bx^{IC}) \\ \mL_s[u_s](\bz; \btheta_s(\tau)) - f(\bz) \end{bmatrix}.
\end{aligned}
\]
with $\bK^s_{11}(\tau) \in \R^{M\times M}$, $\bK_{12}(\tau) \in \R^{M\times N}$ ($\bK^s_{12}(\tau) = (\bK^s_{21}(\tau))^\top$), and $\bK^s_{22}(\tau) \in \R^{N \times N}$ and the entries are given as
\[
\begin{cases}
    (\bK^s_{11})_{i, j}(\tau) &= \inp{\partial_{\btheta}u_s(\bzeta_i; \btheta_s(\tau))\,, \partial_{\btheta}u_s(\bzeta_j; \btheta_s(\tau))} \\[0.2cm]
    (\bK^s_{12})_{i, j}(\tau) &= \inp{\partial_{\btheta}u_s(\bzeta_i; \btheta_s(\tau))\,, \partial_{\btheta}\mL_s[u_s](\bz_j; \btheta_s(\tau))} \\[0.2cm]
    (\bK^s_{22})_{i, j}(\tau) &= \inp{\partial_{\btheta}\mL_s[u_s](\bz_i; \btheta_s(\tau))\,, \partial_{\btheta}\mL_s[u_s](\bz_j; \btheta_s(\tau))} \\
\end{cases}.
\]
Meanwhile, for $u_h$, we have
\[
\frac{\dif \mL_h[u_h](\bz; \btheta_h(\tau))}{\dif\tau}
= -\bK^h(\tau)
\left(
\mL_h[u_h](\bz; \btheta_h(\tau)) - f(\bz) + \psi_t(\bz)
\right),
\quad \bK^h(\tau) \in \R^{N \times N},
\]
% \[
% \frac{\dif \mL_h[u_h](\bz; \btheta_h(\tau))}{\dif\tau} = -\bK^h(\mL_h[u_h](\bz; \btheta_h(\tau)) - f(\bz) - \psi_t(\bz)), \quad \bK^h \in \R^{N \times N},
% \]
with the entries of $\bK^h$ being
\[
(\bK^h)_{i, j}(\tau) = \inp{\partial_{\btheta}\mL_h[u_h](\bz_i; \btheta_h(\tau))\,, \partial_{\btheta}\mL_h[u_h](\bz_j; \btheta_h(\tau))}.
\]
% Setting $\bK^s = \begin{bmatrix} \bK^s_{11} & \bK^s_{12} \\ \bK^s_{21} & \bK^s_{22}\end{bmatrix}$, and assuming $\bK^s(\tau) \approx \bK^s(0) = \bK^s_*$ and $\bK^h(\tau) \approx \bK^h(0) = \bK^h_*$, then the decays satisfy
% \[
% \begin{bmatrix} \frac{\dif u_s(\bzeta; \btheta_s(\tau))}{\dif\tau} \\ \frac{\dif \mL_s[u_s](\bz; \btheta_s(\tau))}{\dif\tau}\end{bmatrix} - \begin{bmatrix} u_0(\bx^{IC}) \\ f(\bz) \end{bmatrix} = -e^{-\bK^s_*\tau}\begin{bmatrix} u_0(\bx^{IC}) \\ f(\bz) \end{bmatrix} = -\bQ_s^\top e^{-\bLambda_s\tau}\bQ_s\begin{bmatrix} u_0(\bx^{IC}) \\ f(\bz) \end{bmatrix},
% \]
% where $\bK^s_* = \bQ_s^\top\bLambda_s\bQ_s$; moreover
% \begin{small}
%     \[
% \frac{\dif \mL_h[u_h](\bz; \btheta_h(\tau))}{\dif\tau} - (f(\bz) + \psi_t(\bz)) = -e^{\bK^h_*\tau}(f(\bz) + \psi_t(\bz)) = -\bQ_h^\top e^{-\bLambda_h\tau}\bQ(f(\bz) + \psi_t(\bz)),
% \]
% \end{small}
% where $\bK^h_* = \bQ_h^\top\bLambda_h\bQ_h$.
Setting
\[
\bK^s =
\begin{bmatrix}
\bK^s_{11} & \bK^s_{12} \\[0.2cm]
\bK^s_{21} & \bK^s_{22}
\end{bmatrix},
\]
and assuming
\[
\bK^s(\tau) \approx \bK^s(0) = \bK^s_*,
\qquad
\bK^h(\tau) \approx \bK^h(0) = \bK^h_*,
\]
the soft-constrained residual satisfies
\[
\begin{aligned}
\begin{bmatrix}
u_s(\bzeta; \btheta_s(\tau)) - u_0(\bx^{IC}) \\
\mL_s[u_s](\bz; \btheta_s(\tau)) - f(\bz)
\end{bmatrix}  =
e^{-\bK^s_*\tau}
\begin{bmatrix}
u_s(\bzeta; \btheta_s(0)) - u_0(\bx^{IC}) \\
\mL_s[u_s](\bz; \btheta_s(0)) - f(\bz)
\end{bmatrix}
 =
\bQ_s^\top e^{-\bLambda_s\tau}\bQ_s
\begin{bmatrix}
u_s(\bzeta; \btheta_s(0)) - u_0(\bx^{IC}) \\
\mL_s[u_s](\bz; \btheta_s(0)) - f(\bz)
\end{bmatrix},
\end{aligned}
\]
where $\bK^s_* = \bQ_s^\top\bLambda_s\bQ_s$. Moreover, the hard-constrained residual satisfies
\[
\begin{aligned}
&\mL_h[u_h](\bz; \btheta_h(\tau)) - f(\bz) + \psi_t(\bz) \\[0.2cm] &=
e^{-\bK^h_*\tau}
\left(
\mL_h[u_h](\bz; \btheta_h(0)) - f(\bz) + \psi_t(\bz)
\right) \\[0.2cm]
& =
\bQ_h^\top e^{-\bLambda_h\tau}\bQ_h
\left(
\mL_h[u_h](\bz; \btheta_h(0)) - f(\bz) + \psi_t(\bz)
\right),
\end{aligned}
\]
where $\bK^h_* = \bQ_h^\top\bLambda_h\bQ_h$.
One can see that, in the soft-constrained formulation, the training decay
depends on the spectral properties of the full coupled block matrix $\bK^s$,
which contains both the initial-condition and PDE-residual components.
In contrast, the training decay of the hard-constrained formulation depends
only on the transformed PDE-residual kernel $\bK^h$, since the initial-condition
residual has been eliminated from the training system. This reduction helps
mitigate the spectral imbalance between competing loss components.

% \begin{theorem}[Well Conditioning]\label{thm:wellcond}
% The transformation is well defined when $\phi$ and $\psi$ is $C^1$ in time and have $K^{th}$ order partial derivatives w.r.t $\bx$.  Moreover, $\phi_t(0, \bx) \neq 0$.  Then $u_{nn}$ will satisfy a new PDE
% \[
% \begin{cases}
%     &\phi\partial_tu_{nn} + \mathcal{P}[\psi + \phi u_{nn}] + \phi_tu_{nn} = f - \psi_t, \quad (t, \bx) \in (0, T]\times\Omega, \\
%     &u_{nn}(0, \bx) = \frac{u_t(0, \bx) - \psi_t(0, \bx)}{\phi_t(0, \bx)}, \quad \bx \in \bar\Omega.
% \end{cases}
% \]
% \end{theorem}

\begin{theorem}[Hessian Characterization and Bounds]\label{thm:hessian}
Assume that $\mP$ is linear and that $\phi(t,\bx)=\phi(t)\neq 0$
for $t>0$. Consider a linear-in-parameters representation of the
network output, or equivalently its NTK linearization about a fixed
parameter initialization,
\[
    u(\by;\btheta)=\h^\top(\by)\btheta,
    \qquad
    \h(\by)=
    \begin{bmatrix}
        h_1(\by) & \cdots & h_D(\by)
    \end{bmatrix}^{\top},
\]
% Assume that $\mP$ is linear and $\phi(t, \bx) = \phi(t) \neq 0$ for $t > 0$, then for a one-hidden layer network, we have for a network $u(\by; \theta)$ with $\btheta \in \R^{D}$ being the network parameters, the following relationship is satisfied
% \[
% u(\by; \btheta) = \h^\top(\by)\btheta, \quad \h = \begin{bmatrix} h_1(\by) & \cdots & \h_D(\by) \end{bmatrix}^\top,
% \]
where $\by = \bz = (t, \bx) \in (0, T]\times\Omega$ or $\by = \bzeta = (0, \bx^{IC}) \in \{0\}\times\bar\Omega$.  Therefore the two losses become
\[
\begin{aligned}
\text{Loss}_s(u(\by;\btheta)) &= \frac{1}{2N}\norm{\bA\btheta - \f}_{\ell^2(\R^{N})}^2 + \frac{1}{2M}\norm{\bC\btheta - \bu_0}_{\ell^2(\R^{M})}^2, \; \text{(soft constrained)}\\
\text{Loss}_h(u(\bz;\btheta)) &= \frac{1}{2N}\norm{\bLambda_{\phi}\bA\btheta + \bLambda_{\phi_t}\bB\btheta - \tilde\f}_{\ell^2(\R^{N})}^2, \; \text{(hard constrained)}
\end{aligned}
\]
where
\[
\bA = \begin{bmatrix}\mL_s[\h^\top(\bz_1)] \\ \vdots \\ \mL_s[\h^\top(\bz_N)]\end{bmatrix}, \quad \bB = \begin{bmatrix}\h^\top(\bz_1) \\ \vdots \\ \h^\top(\bz_N)\end{bmatrix}, \quad \text{and} \quad \bC = \begin{bmatrix}\h^\top(\bzeta_1) \\ \vdots \\ \h^\top(\bzeta_M)\end{bmatrix},
\]
where $\bx^{IC}_i \in \bar\Omega$; and for the diagonal matrices $\bLambda_{\phi}$ and $\bLambda_{\phi_t}$
\[
\bLambda_{\phi} = \text{diag}(\phi(\bz_1), \cdots, \phi(\bz_N)) \quad \text{and} \quad \bLambda_{\phi_t} = \text{diag}(\phi_t(\bz_1), \cdots, \phi_t(\bz_N)),
\]
and the vectors
\[
\f = \begin{bmatrix}f(\bz_1) \\ \vdots \\ f(\bz_N) \end{bmatrix}, \quad \tilde\f = \begin{bmatrix}f(\bz_1) - \mL_s[\psi](\bz_1) \\ \vdots \\ f(\bz_N) - \mL_s[\psi](\bz_N) \end{bmatrix}, \quad \text{and} \quad \bu_0 = \begin{bmatrix}u_0(\bx_1^{IC}) \\ \vdots \\ u_0(\bx_M^{IC}) \end{bmatrix}.
\]
Therefore the Hessians of the losses are
\[
\begin{aligned}
&H_s(\btheta) = \frac{1}{N}\bA^\top\bA + \frac{1}{M}\bC^\top\bC,\;
H_h(\btheta)\\&= \frac{1}{N}
\Bigl(
\bA^\top\bLambda_{\phi}^2\bA
+
\bA^\top\bLambda_{\phi}\bLambda_{\phi_t}\bB
+
\bB^\top\bLambda_{\phi_t}\bLambda_{\phi}\bA
+
\bB^\top\bLambda_{\phi_t}^2\bB
\Bigr),\\
\end{aligned}
\]
% \[
% \begin{aligned}
% H_h(\btheta)
% &=
% \frac{1}{N}
% \left(
% \bLambda_{\phi}\bA
% +
% \bLambda_{\phi_t}\bB
% \right)^\top
% \left(
% \bLambda_{\phi}\bA
% +
% \bLambda_{\phi_t}\bB
% \right)=
% \frac{1}{N}
% \Bigl(
% \bA^\top\bLambda_{\phi}^2\bA
% +
% \bA^\top\bLambda_{\phi}\bLambda_{\phi_t}\bB
% +
% \bB^\top\bLambda_{\phi_t}\bLambda_{\phi}\bA
% +
% \bB^\top\bLambda_{\phi_t}^2\bB
% \Bigr).
% \end{aligned}
% \]

Then the Hessians are bounded as
\[
\begin{aligned}
    \norm{H_s}_2 \le \frac{\norm{\mL_s}_2^2}{N}\norm{\bB}_2^2 + \frac{1}{M}\norm{\bC}_2^2, \quad
    \norm{H_h}_2 \le \frac{1}{N}(C_1\norm{\mL_s}_2^2 + 2C_2\norm{\mL_s}_2 + C_3)\norm{\bB}_2^2.
\end{aligned}
\]
where $\norm{\mL_s}_2$ denotes an induced discrete operator bound
on the sampled linearized trial space such that
\[
    \norm{\bA}_2 \leq \norm{\mL_s}_2\norm{\bB}_2,
\]
and
\[
C_1 = \max_{\bz}\phi^2(\bz), \quad C_2 = \max_{\bz}|\phi(\bz)\phi_t(\bz)|, \quad \text{and} \quad C_3 = \max_{\bz}\phi_t^2(\bz).
\]
\end{theorem}
% \begin{remark}
% Notice the competition between the PDE residual points and the initial conditions points (reflected in $\bB$ and $\bC$) in the bound for $H_s$, where for $H_h$, it is all about the operator $\mL_s$ and the scaling by $\phi$ and $\phi_t$.  The bounds will become much more complex when $\mP$ is not linear or $\phi$ also depends on $\bx$.

% % In $\text{Loss}_s$, PDE and IC residuals compete, producing anisotropic curvature (two objectives).
% % In $\text{Loss}_h$, all residuals are expressed through a single PDE operator scaled by $\phi,\phi_t$,
% % yielding a smoother, better-conditioned optimization landscape.

% \end{remark}

\begin{remark} The soft-constrained Hessian contains separate contributions from the PDE residual and the initial-condition mismatch, represented by $\bA^\top\bA$ and $\bC^\top\bC$, respectively. In contrast, the hard-constrained Hessian is determined entirely by the transformed PDE residual through $\bLambda_{\phi}\bA+\bLambda_{\phi_t}\bB$, since the initial-condition residual has been eliminated from the loss. This structural simplification is consistent with the NTK decay analysis above. When $\mP$ is nonlinear or $\phi$ also depends on $\bx$, additional terms arise and the corresponding Hessian characterization becomes more involved. \end{remark}

\section{Examples}\label{section:example}

In this section, we examine the hard-constrained PINN (HC-PINN) on a diverse set of prototypical stiff, time-dependent PDEs with periodic boundary conditions, aiming to demonstrate its robustness in addressing a broad class of challenging problems.   The
benchmarks in Section~\ref{sec:ac} -- \ref{sec:la} all use periodic boundary conditions.
Section~\ref{sec:dirichlet} additionally verifies that the hard-IC
mechanism extends to non-periodic settings by reporting a Dirichlet
Allen--Cahn run with a standard MLP (no Fourier features). These equations often involve strong nonlinearities, stiffness, and parameter sensitivities, which pose significant challenges for data-driven methods such as PINNs. Nonetheless, our results illustrate that the proposed framework can effectively overcome these difficulties and accurately approximate solutions across these regimes.

We present the results in five major categories: $(I)$ an ablation study across $7$ stiff PDEs with different sources of stiffness, including sharp phase transitions in Allen--Cahn, higher-order derivatives (fourth order) in Cahn--Hilliard and Kuramoto--Sivashinsky, coupled outputs in Gray--Scott and Schr\"odinger systems, and high-frequency modes in fast-moving linear transport; $(II)$ comparisons with state-of-the-art PINNs, including the Causal PINN~\cite{wang2022respecting} and Enhanced RBA PINN~\cite{anagnostopoulos2024residual}, on Allen--Cahn, Cahn--Hilliard, and Kuramoto--Sivashinsky; $(III)$ application of HC-PINNs to a $2D$ Allen--Cahn equation; $(IV)$ sensitivity analysis with respect to the parameter $m$ in the Allen--Cahn equation under two different initial conditions; and $(V)$, as an additional robustness check beyond the periodic Fourier-feature setting, a non-periodic Dirichlet Allen--Cahn run with a plain MLP, reported in Section~\ref{sec:dirichlet}.

\subsection{Implementation details}\label{section:B}
All experiments were implemented in PyTorch ($v2.2$) and executed on CUDA-enabled NVIDIA GPUs as well as Apple Silicon devices via the MPS backend. Our code runs efficiently across Linux (Ubuntu $22.04$), macOS (Sonoma), and Windows $11$ platforms. Primary training was conducted on an NVIDIA $A100$ GPU ($40$ GB memory) and a MacBook Pro with an Apple M1 Pro chip ($16$-core GPU, $32$ GB unified memory).

For each benchmark, we used a fully connected neural network with $7$ hidden layers, each containing $32$ neurons and tanh activation functions, to represent the latent solution. The same initial weights were used across all ablation experiments for consistency.

Model training was carried out using a two-stage strategy: we first trained with the Adam optimizer using an initial learning rate of $10^{-3}$ for $50{,}000$ steps, followed by fine-tuning with the L-BFGS-B optimizer with a learning rate of $1$ for $5{,}000$ steps. No hyperparameter sweeps were performed, and the total number of training iterations was fixed across all experiments.

This work does not rely on external datasets. We used Latin Hypercube Sampling (LHS) to generate a fixed set of $16{,}384$ interior collocation points, along with $128$ initial data points and $128$ boundary data points. This sampling configuration was kept consistent across all experiments to ensure comparability. For PDE systems, we employed a mini-batching strategy with $5$ batches of size $3{,}280$ to improve memory efficiency during training.

To facilitate training and better incorporate the initial condition into the solution structure, we adopt the transformation
% \[
% \tilde{u}(t, \bx) = e^{-Ct} u_0(\bx) + t\,u_{nn}(t, \bx),
% \quad \psi = e^{-Ct} u_0(\bx), \quad \phi = t, \quad (t, \bx) \in [0, T]\times\bar\Omega,
% \]
\[
\tilde{u}(t, \bx) = e^{-Ct} u_0(\bx)
+ (t/T)\,u_{nn}(t, \bx),
\quad \psi = e^{-Ct} u_0(\bx),\]\[\hspace{-3.3cm}
\quad \phi(t) = t/T,
\quad (t, \bx) \in [0, T]\times\bar\Omega,
\]
where $u_0(\bx)$ denotes the initial condition. The
factor $\phi(t)=t/T$ is equivalent to training the network on the
normalized time coordinate $t/T\in[0,1]$, which is what we use in
practice. This setting keeps $\phi$ bounded on the simulation interval for any
final time $T$, while preserving the non-degeneracy
$\phi_t(0,x)=1/T\neq 0$ required by Theorem~\ref{thm:wellposed}.
We test two values of $C \in \{0.1, 1\}$ to assess how the decay rate affects the learned dynamics. The transformed PINN $u_{nn}$ is trained by minimizing the PDE-residual loss,
while the initial condition is enforced exactly by construction now.
% The transformed PINN $u_{nn}$ is trained to minimize a loss function that accounts for the PDE residuals, which now implicitly enforces the IC/BC.

% Ablation studies were conducted by disabling one component of the method at a time while keeping all other settings fixed, to evaluate its contribution using relative $L^2$ and relative $L^1$ error metrics. These measure discrepancies between predicted and reference solutions in average (energy-based) and absolute (mass-based) senses, respectively, and are normalized to enable fair comparisons across different scales.

Ablation studies were conducted by disabling one component of the method at a time
while keeping all other settings fixed, to evaluate its contribution using relative
$L^2$ and relative $L^1$ error metrics. These quantify quadratic/energy-error and
absolute-error discrepancies, respectively, and are normalized to enable fair
comparisons across different scales.

We also compare HC-PINN with causal training (MLP) and enhanced RBA using
the same data and a total of $50{,}000$ iterations, instead of the
$300{,}000$ iterations used in the original works. This
choice is made deliberately to enforce an equal compute budget across all
methods, so that differences in accuracy reflect the training schemes
themselves rather than disparities in optimization effort; we acknowledge
that this reduced budget may lower the performance of the baseline methods
relative to their originally reported values, and we report those original
$300{,}000$-iteration values explicitly in Table~\ref{tab:comparison} for
transparency. For all baselines we use the recommended/default
hyperparameters from their respective original works, and, consistent with
the no-sweep protocol stated above, no additional tuning was performed for
either the baselines or HC-PINN. For the sensitivity study, experiments with different $m$ values are conducted on the same set of points, using the same number of iterations and architecture. For $2D$ training, we use $4{,}096$ initial collocation points, $10{,}000$ residual collocation points per segment, and $800$ boundary points per segment. The architecture remains the same as in the $1D$ setting.

\subsection{Allen--Cahn PDE}\label{sec:ac}
The Allen--Cahn type PDE is a classical phase-field model that has been widely used to study phase separation phenomena \cite{bazant2017thermodynamic}. It has numerous practical applications across a range of fields, including materials science \cite{allen1979microscopic,shen2010numerical}, biological systems \cite{hyman2014liquid,takatori2015towards}, and electrochemical systems \cite{horstmann2013precipitation,tian2015electrochemically}.
We begin with the $1$-dim Allen--Cahn equation with periodic boundary conditions formulated as follows:
\begin{equation}\label{eq:ac}
\begin{aligned}
&u_t - \gamma_1 u_{xx} +\gamma_2(u^3 - u) = 0,  \, (t, x) \in (0, T]\times(a, b), \\
&u(0, x) = u_0(x), \, x \in [a, b], \quad u(t,a) = u(t,b), \, t \in [0, T],
\end{aligned}
\end{equation}
 where $\gamma_1,\gamma_2 > 0, T > 0, a < b$ are prescribed constants. As $\gamma_2$ increases, the transition interface of the solutions is sharper, which makes it harder to solve the AC equation numerically.

\subsubsection{Case \texorpdfstring{$I$}{I}}
We first tested on Allen--Cahn equation with the following, $u_0(x) = x^2\cos(\pi x)$, $T = 1$, $a = -1$, $b = 1$, $\gamma_1 = 0.001$, and $\gamma_2 =5$.

% \begin{figure}[h]
% \centering
% \includegraphics[width=0.48\textwidth]{figs/fig6_AC.png}
% \caption{Exact solution of  \textit{Allen-Cahn}  with the corresponding best network prediction and the absolute error difference.
% % The resulting relative $L^2$ error is $1.18\times10^{-3}$ and the relative $L^1$ error is $5.34\times 10^{-4}$. The hyper-parameter configuration is specified in Section \ref{section:3}.
% }
% \label{fig:ac_utx}
% \end{figure}

% \begin{figure}[h]
% \centering
% \includegraphics[width=0.48\textwidth]{original_doc/fig6_AC_t.png}
% \caption{Best solution: \textit{Allen-Cahn}.
% % Best solutions of the \textit{Allen-Cahn Equation (Case I)} $\gamma_2=5$. The resulting relative $L^2$ error is $1.18\times10^{-3}$ and the relative $L^1$ error is $5.34\times 10^{-4}$. The hyper-parameter configuration is specified in Section \ref{section:3}.
% }
% \label{fig:ac}
% \end{figure}
\begin{table}[h]
%title of the table
\setlength{\tabcolsep}{0.75mm}
\centering
\resizebox{\textwidth}{!}{
\renewcommand{\arraystretch}{1.5}
\begin{tabular}{cccccccc}
\toprule
\multicolumn{4}{c}{\textbf{Ablation Settings}}&\multicolumn{2}{c}{\textbf{Case $I$} }&\multicolumn{2}{c}{\textbf{Case $II$} }\\ [0.5ex]
\cmidrule(lr){1-4} \cmidrule(lr){5-6} \cmidrule(lr){7-8}
\textbf{PBC} & \textbf{IC1} & \textbf{IC2} & \textbf{SA}& \textbf{Rel. $L^2$ error}& \textbf{Rel. $L^1$ error}& \textbf{Rel. $L^2$ error}& \textbf{Rel. $L^1$ error}\\[0.2ex]
\hline
 % Entering row contents Midnight&7&-3&
\Checkmark & \Checkmark & \XSolidBrush &\Checkmark & $8.29\times 10^{-4}$ & $5.34\times 10^{-4}$ & $3.53\times 10^{-4}$ & $2.14\times 10^{-4}$\\[0.5ex]
\Checkmark & \XSolidBrush & \Checkmark &\Checkmark & $7.64\times 10^{-2}$ & $1.66\times 10^{-2}$ & $7.04\times 10^{-4}$ & $3.26\times 10^{-4}$\\[0.5ex]
\Checkmark & \Checkmark & \XSolidBrush &\XSolidBrush & $1.24\times 10^{-3}$ & $6.10\times 10^{-4}$ &$8.32\times 10^{-4}$ & $3.72\times 10^{-4}$\\[0.5ex]
\Checkmark & \XSolidBrush & \Checkmark &\XSolidBrush & $2.86\times 10^{-2}$ & $8.20\times 10^{-3}$ & $9.80\times 10^{-4}$ & $3.91\times 10^{-4}$\\[0.5ex]
\Checkmark& \XSolidBrush  & \XSolidBrush &\XSolidBrush & $9.99\times 10^{-1}$ & $9.99\times 10^{-1}$ & $9.99\times 10^{-1}$ & $9.99\times 10^{-1}$\\[0.5ex]
\Checkmark& \XSolidBrush  & \XSolidBrush &\Checkmark& $9.99\times 10^{-1}$ & $9.99\times 10^{-1}$ & $9.99\times 10^{-1}$ & $9.99\times 10^{-1}$\\[0.5ex]
\XSolidBrush  & \Checkmark & \XSolidBrush  &\XSolidBrush & $1.05\times 10^{-2}$ & $2.65\times 10^{-3}$ & $5.11\times 10^{-2}$ & $1.21\times 10^{-2}$\\[0.5ex]
\XSolidBrush  & \Checkmark & \XSolidBrush  &\Checkmark & $1.86\times 10^{-3}$ & $7.31\times 10^{-4}$ & $5.03\times 10^{-2}$ & $8.27\times 10^{-3}$\\[0.5ex]
\XSolidBrush  &   \XSolidBrush&\Checkmark  &\XSolidBrush & $1.88\times 10^{-2}$ & $4.71\times 10^{-3}$  & $4.58\times 10^{-2}$ & $1.10\times 10^{-2}$\\[0.5ex]
\XSolidBrush & \XSolidBrush   & \Checkmark &\Checkmark & $4.34\times 10^{-3}$ & $1.24\times 10^{-3}$  & $3.03\times 10^{-2}$ & $4.25\times 10^{-3}$\\[0.5ex]
\XSolidBrush  &  \XSolidBrush  & \XSolidBrush  &\Checkmark & $5.12\times 10^{-1}$ & $3.18\times 10^{-1}$ & $4.02\times 10^{-1}$ & $1.42\times 10^{-1}$\\[0.5ex]
\XSolidBrush  &  \XSolidBrush  & \XSolidBrush  &\XSolidBrush & $9.99\times 10^{-1}$ & $9.98\times 10^{-1}$ & $5.11\times 10^{-1}$ & $3.28\times 10^{-1}$\\[0.5ex]
\bottomrule
\end{tabular}}
\vspace{0.2cm}
\caption{\textit{Allen--Cahn (Case $I$ and Case $II$):} Relative $L^2$ and $L^1$ errors for an ablation study illustrating the impact of disabling individual components of the proposed technique and training pipeline.}
\label{tab:ac}
\end{table}
For this example, we first test the baseline PINN approach \cite{raissi2019physics}. As shown in Table~\ref{tab:ac}, the baseline PINN alone fails to accurately solve the Allen--Cahn equation, motivating the development of a more effective training pipeline. To address this, we employ all proposed techniques (except mini-batching) and conduct an ablation study by disabling individual components one at a time. Additionally, we introduce two different decay functions for the initial condition, $e^{-1t}$ and $e^{-0.1t}$, to investigate their influence on the solution. The results are summarized in Table \ref{tab:ac}. The full scheme, including self-adaptive loss balancing, achieves the best performance, reaching a relative $L^2$ error of $8.29\times 10^{-4}$ and relative $L^1$ error of $5.34\times 10^{-4}$. In particular, using the faster-decaying initial condition ($IC1$) leads to better results than the slower-decaying case ($IC2$), suggesting that the initial condition plays a more critical role during the early stages of training but becomes less dominant over time. The ablation results further highlight the importance of each methodological component. Removing any single part generally leads to performance degradation, with the most severe impact observed when the initial condition is omitted, resulting in a relative error close to $0.99$ in both $L^1$ and $L^2$ norms.
Figure~\ref{fig:ac} presents cross-sectional comparisons at different times, further confirming the excellent agreement between the predicted and true solutions across the entire time horizon.

\subsubsection{Case \texorpdfstring{$II$}{II}}
For the second case, we change the IC to have highly oscillatory data by $u(0, x) = x^2\sin(2\pi x)$. Other parameters are $T = 1$, $a = -1$, $b = 1$, $\gamma_1 = 0.001$ and $\gamma_2=4$.

\begin{figure}[ht]
\centering
\includegraphics[width=0.75\textwidth]{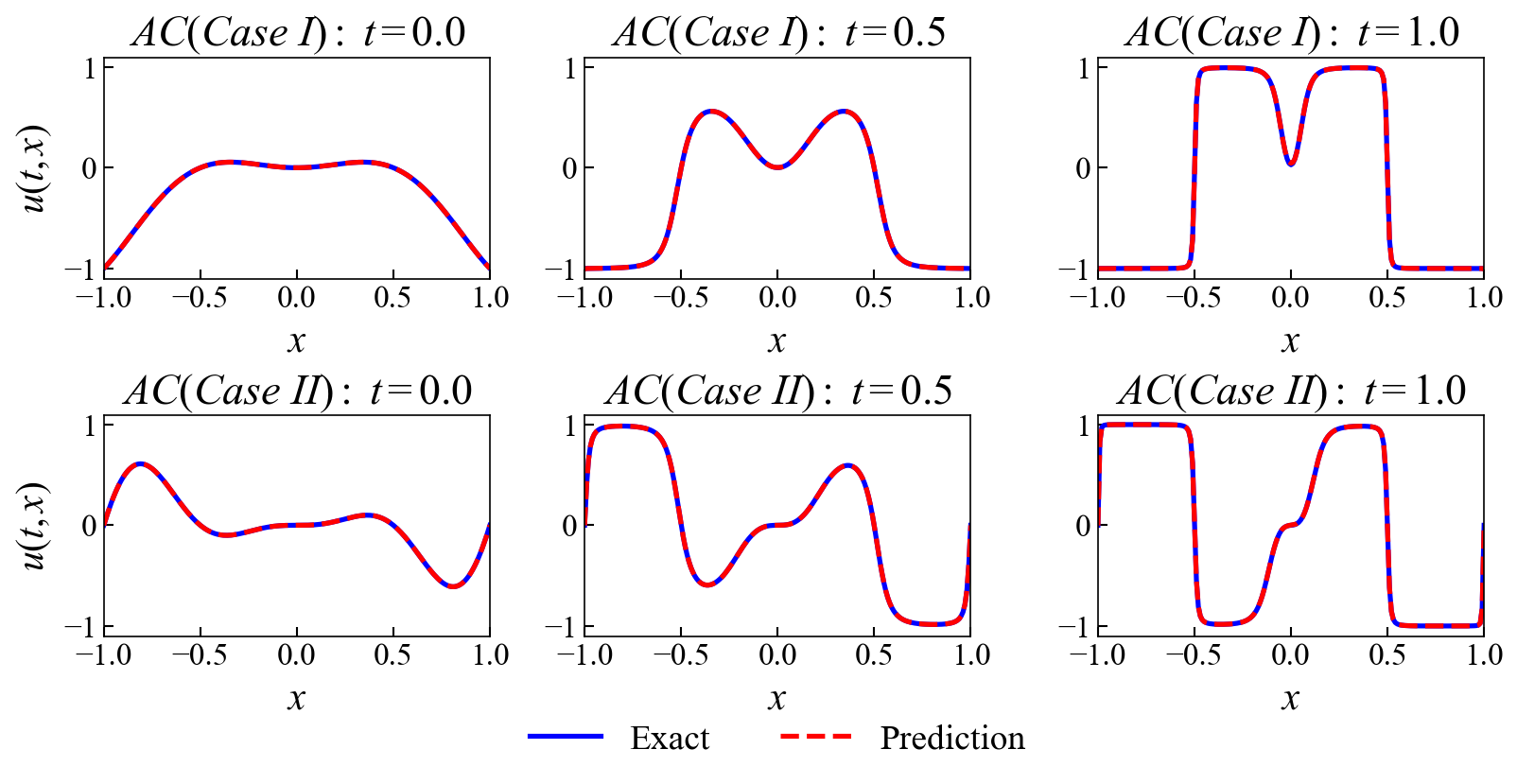}
\caption{
Top: Best solution of the \textit{Allen--Cahn Equation (Case $I$)} $\gamma_2=5$,
obtained with the full scheme (PBC, IC1, SA). The resulting
relative $L^2$ error is $8.29\times10^{-4}$ and the relative
$L^1$ error is $5.34\times 10^{-4}$.
Bottom: Best solution of the \textit{Allen--Cahn Equation (Case $II$)}
$\gamma_2=4$, obtained with the full scheme (PBC, IC1, SA).
The resulting relative $L^2$ error is $3.53\times10^{-4}$ and the relative
$L^1$ error is $2.14\times 10^{-4}$.
Full $(x,t)$ network predictions are in Fig.~\ref{fig:all}.
}
\label{fig:ac}
\end{figure}

We first train the baseline PINN model, which produces a large relative $L^2$ error of $0.511$ and relative $L^1$ error of $0.328$, demonstrating its inadequacy in resolving the sharp interface dynamics (Table~\ref{tab:ac}). However, once the initial condition enforcement and periodic boundary neural networks are incorporated, the performance improves drastically, achieving a relative $L^2$ error of $8.32\times 10^{-4}$ and relative $L^1$ error of $3.72\times 10^{-4}$.  Adding self-adaptive
weighting further improves the full-scheme errors to $3.53\times 10^{-4}$
and $2.14\times 10^{-4}$ (Figure~\ref{fig:ac}, bottom). Further gains are realized when using the self-adaptive loss weighting strategy, although the improvement is marginal compared to the dominant contribution from enforcing the initial and boundary conditions. Interestingly, unlike the previous case, there is minimal performance difference between the schemes using $IC1$ and $IC2$, suggesting that the impact of initial condition decay becomes less significant when strong oscillations are already present at the beginning.

\subsection{Kuramoto--Sivashinsky PDE}\label{sec:ks}
\begin{figure}[h]
\centering
\includegraphics[width=0.7\textwidth]{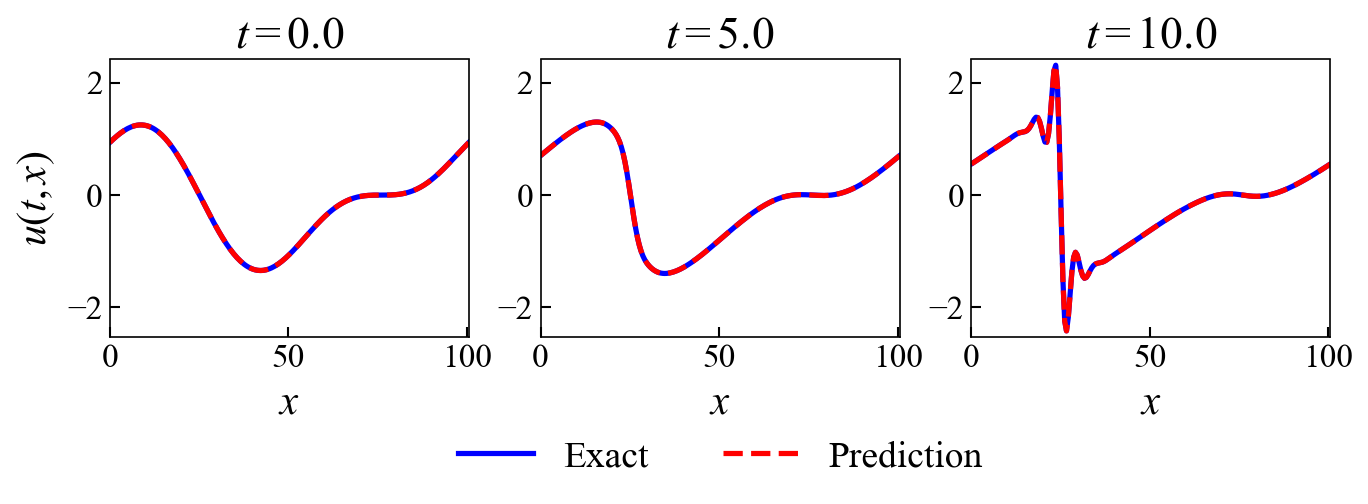}
\caption{Best solution: \textit{Kuramoto--Sivashinsky}, obtained
with the configuration (PBC, IC2, SA). The resulting relative $L^2$ error is
$5.37\times10^{-4}$ and the relative $L^1$ error is
$1.02\times 10^{-3}$. Full $(x,t)$ network predictions are in
Fig.~\ref{fig:all}.}
\label{fig:ks}
\end{figure}

Kuramoto--Sivashinsky (KS) equation is a fourth-order nonlinear partial differential equation, known for its chaotic behavior \cite{kudryashov1990exact,michelson1986steady}.  It has served as a pivotal model in the study of complex dynamical phenomena observed in various physical systems. Consequently, the KS equation provides an excellent case study to demonstrate how PINNs can effectively simulate chaotic dynamics. The equation is given by:
 \begin{equation}\label{eq:ks}
\begin{aligned}
u_t& =  -u_{xx} - u_{xxxx} - uu_x,  \quad (t, x) \in (0, T]\times (a, b), \\
u(0, x) &= \cos(\frac{x}{16})(1 + \sin\frac{x-1}{16}), \quad x \in [a, b], \quad u(t,a) = u(t,b), \, t \in [0, T],
\end{aligned}
\end{equation}
 where $T=10, a=0, b=32\pi$.  The best results are shown in Figure~\ref{fig:ks}.

\begin{table}[t]
%title of the table
\centering
\renewcommand{\arraystretch}{1.4}
\begin{tabular}{cccccc}
\toprule
\multicolumn{4}{c}{\textbf{Ablation Settings}}&\multicolumn{2}{c}{\textbf{Performance} }\\ [0.5ex]
\cmidrule(lr){1-4} \cmidrule(lr){5-6}
\textbf{PBC} & \textbf{IC1} & \textbf{IC2} & \textbf{SA}& \textbf{Rel. $L^2$ error}& \textbf{Rel. $L^1$ error}\\[0.2ex]
\hline
 % Entering row contents Midnight&7&-3&
\Checkmark & \Checkmark & \XSolidBrush &\Checkmark & $4.04\times 10^{-2}$ & $1.46\times 10^{-2}$\\[0.5ex]
\Checkmark & \XSolidBrush & \Checkmark &\Checkmark & $5.37\times 10^{-4}$ & $1.02\times 10^{-3}$ \\[0.5ex]
\Checkmark & \Checkmark & \XSolidBrush &\XSolidBrush & $1.68\times 10^{-2}$ & $8.94\times 10^{-3}$\\[0.5ex]
\Checkmark & \XSolidBrush & \Checkmark &\XSolidBrush & $1.30\times 10^{-2}$ & $7.79\times 10^{-3}$\\[0.5ex]
\Checkmark& \XSolidBrush  & \XSolidBrush &\XSolidBrush & $9.84\times 10^{-3}$ & $6.36\times 10^{-3}$\\[0.5ex]
\Checkmark& \XSolidBrush  & \XSolidBrush &\Checkmark& $2.50\times 10^{-3}$ & $1.81\times 10^{-3}$\\[0.5ex]
\XSolidBrush  & \Checkmark & \XSolidBrush  &\XSolidBrush & $2.13\times 10^{-1}$ & $9.74\times 10^{-2}$\\[0.5ex]
\XSolidBrush  & \Checkmark & \XSolidBrush  &\Checkmark & $2.00\times 10^{-1}$ & $1.04\times 10^{-1}$\\[0.5ex]
\XSolidBrush  &   \XSolidBrush&\Checkmark  &\XSolidBrush & $1.54\times 10^{-2}$ & $1.21\times 10^{-2}$\\[0.5ex]
\XSolidBrush & \XSolidBrush   & \Checkmark &\Checkmark & $9.60\times 10^{-3}$ & $8.75\times 10^{-3}$\\[0.5ex]
\XSolidBrush  &  \XSolidBrush  & \XSolidBrush  &\Checkmark & $2.11\times 10^{-1}$ & $1.06\times 10^{-1}$\\[0.5ex]
\XSolidBrush  &  \XSolidBrush  & \XSolidBrush  &\XSolidBrush & $3.01\times 10^{-1}$ & $1.97\times 10^{-1}$\\[0.5ex]
\bottomrule
\end{tabular}
\vspace{0.4cm}
\caption{\textit{Kuramoto--Sivashinsky:} Relative $L^2$ and $L^1$ errors for an ablation study illustrating the impact of disabling individual components of the proposed technique and training pipeline.} \label{tab:ks}
\end{table}

As in Table~\ref{tab:ks}, solving this PDE system using a baseline PINN alone proves to be exceptionally challenging. The incorporation of periodic neural networks for boundary conditions, paired with initial condition \textit{IC2}, consistently yields superior results compared to \textit{IC1}. This suggests that, in contrast to the Allen--Cahn equation discussed in Section \ref{sec:ac}, the choice of initial condition plays a more significant role in the training process for this more complex system. This observation underscores the importance of carefully selecting both the boundary conditions and initial conditions when training neural networks for chaotic or higher-order PDE systems, where these factors contribute notably to the accuracy and stability of the solution.
\subsection{Cahn--Hilliard equation}
The Cahn--Hilliard equation, a variant of the Allen--Cahn equation, also models phase separation. This process describes how the components of a binary fluid spontaneously separate into distinct domains, each composed of a single component \cite{cahn1958free, chan2025phase, kim2016basic, lee2014physical, miranville2017cahn}. The equation is expressed as:
\begin{equation}
    \begin{aligned}
u_t &= \epsilon_1(-u_{xx}  - \epsilon_2u_{xxxx} + (u^3)_{xx}), \quad(t, x) \in (0, T]\times(-L, L),\\
u(0, x)& = u_0(x), \quad  x \in [-L, L], \\
u(t, -L) &= u(t, L), \quad  t \in [0, T],
\end{aligned}
\end{equation}
where $\epsilon_1=10^{-2}, \epsilon_2=10^{-4}, T=1, L=1$. We specify the initial condition as $u_0(x) = -\cos(2\pi x)$. This PDE involves higher-order derivatives, making it more challenging to solve compared to the Allen--Cahn equation.
%
%  \begin{figure}[h]
% \centering
% \includegraphics[width=0.6\textwidth]{original_doc/fig3_CH_BL+PBC+IC2+SA.png}
% \caption{
% Exact solution of the \textit{Cahn-Hilliard} with the corresponding best network prediction and the absolute error difference.
% % The resulting relative $L^2$ error is $9.75\times10^{-4}$ and the relative $L^1$ error is $5.65\times 10^{-4}$. The hyper-parameter configuration is specified in Section \ref{section:3}.
% }
% \label{fig:ch_utx}
% \end{figure}

% \begin{figure}[h]
% \centering
% \includegraphics[width=0.55\textwidth]{fig_CH_BL+PBC+IC2+SA.png}
% \caption{Best solution of the \textit{Cahn-Hilliard}. The resulting relative $L^2$ error is \textcolor{red}{$9.75\times10^{-4}$} and the relative $L^1$ error is \textcolor{red}{$5.65\times 10^{-4}$}. Full $(x,t)$  network predictions are in Fig.~\ref{fig:all}.
% }
% \label{fig:ch}
% \end{figure}

\begin{figure}[h]
\centering
\includegraphics[width=0.65\textwidth]{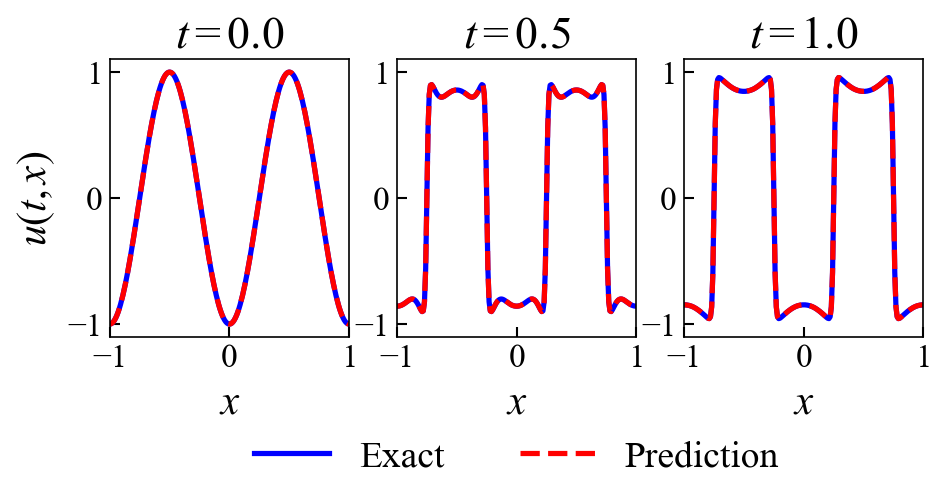}
\caption{Best solution of the \textit{Cahn--Hilliard}, obtained
with the configuration (PBC, IC2, SA). The resulting relative $L^2$ error is
$9.75\times10^{-4}$ and the relative $L^1$ error is
$5.65\times 10^{-4}$. Full $(x,t)$ network predictions are
in Fig.~\ref{fig:all}.}
\label{fig:ch}
\end{figure}

In this experiment, the full scheme yields the best results. Additionally, methods using initial condition \textit{IC2} outperform those with \textit{IC1}. Except for the full schemes and the experiment with periodic neural networks for boundary conditions and initial condition \textit{IC2}, none of the other approaches produce accurate solutions. The errors for all experiments are summarized in Table \ref{tab:ch}. The best results are illustrated in Figure~\ref{fig:ch}.
\begin{table}[h]
%title of the table
\centering
\renewcommand{\arraystretch}{1.5}
\begin{tabular}{cccccc}
\toprule
\multicolumn{4}{c}{\textbf{Ablation Settings}}&\multicolumn{2}{c}{\textbf{Performance} }\\ [0.5ex]
\cmidrule(lr){1-4} \cmidrule(lr){5-6}

\textbf{PBC} & \textbf{IC1} & \textbf{IC2} & \textbf{SA}& \textbf{Rel. $L^2$ error}& \textbf{Rel. $L^1$ error}\\[0.2ex]
\hline

 % Entering row contents Midnight&7&-3&
\Checkmark & \Checkmark & \XSolidBrush &\Checkmark & $1.06\times 10^{-3}$ & $5.61\times 10^{-4}$\\[0.5ex]
\Checkmark & \XSolidBrush & \Checkmark &\Checkmark & $9.75\times 10^{-4}$ & $5.65\times 10^{-4}$ \\[0.5ex]
\Checkmark & \Checkmark & \XSolidBrush &\XSolidBrush & $2.50\times 10^{-1}$ & $5.87\times 10^{-2}$\\[0.5ex]
\Checkmark & \XSolidBrush & \Checkmark &\XSolidBrush & $3.87\times 10^{-3}$ & $1.91\times 10^{-3}$\\[0.5ex]
\Checkmark& \XSolidBrush  & \XSolidBrush &\XSolidBrush & $9.99\times 10^{-1}$ & $9.99\times 10^{-1}$\\[0.5ex]
\Checkmark& \XSolidBrush  & \XSolidBrush &\Checkmark& $9.99\times 10^{-1}$ & $9.99\times 10^{-1}$\\[0.5ex]
\XSolidBrush  & \Checkmark & \XSolidBrush  &\XSolidBrush & $3.48\times 10^{-1}$ & $2.94\times 10^{-1}$\\[0.5ex]
\XSolidBrush  & \Checkmark & \XSolidBrush  &\Checkmark & $2.37\times 10^{-1}$ & $1.09\times 10^{-1}$\\[0.5ex]
\XSolidBrush  &   \XSolidBrush&\Checkmark  &\XSolidBrush & $4.07\times 10^{-1}$ & $3.45\times 10^{-1}$\\[0.5ex]
\XSolidBrush & \XSolidBrush   & \Checkmark &\Checkmark & $1.55\times 10^{-2}$ & $8.20\times 10^{-3}$\\[0.5ex]
\XSolidBrush  &  \XSolidBrush  & \XSolidBrush  &\Checkmark & $9.99\times 10^{-1}$ & $9.99\times 10^{-1}$\\[0.5ex]
\XSolidBrush  &  \XSolidBrush  & \XSolidBrush  &\XSolidBrush & $9.99\times 10^{-1}$ & $9.98\times 10^{-1}$\\[0.5ex]
\bottomrule
\end{tabular}
\vspace{0.4cm}
\caption{\textit{Cahn--Hilliard Equation:} Relative $L^2$ and $L^1$ errors for an ablation study illustrating the impact of disabling individual components of the proposed technique and training pipeline.} \label{tab:ch}
\end{table}
\subsection{Gray--Scott equation}
\begin{figure}[h]
    \centering
    % First figure
    \begin{minipage}{0.48\textwidth}
        \centering
        \includegraphics[width=\linewidth]{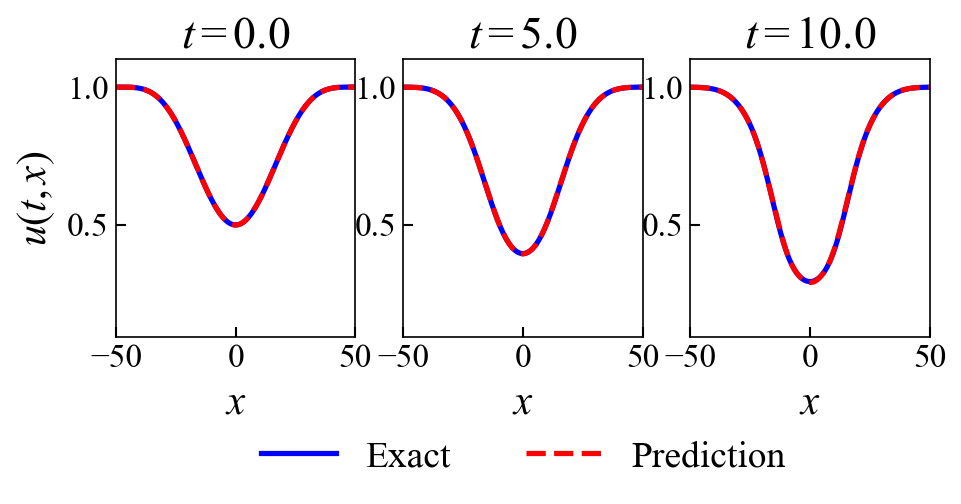}
    \end{minipage}
    \hfill
    % Second figure
    \begin{minipage}{0.48\textwidth}
        \centering
        \includegraphics[width=\linewidth]{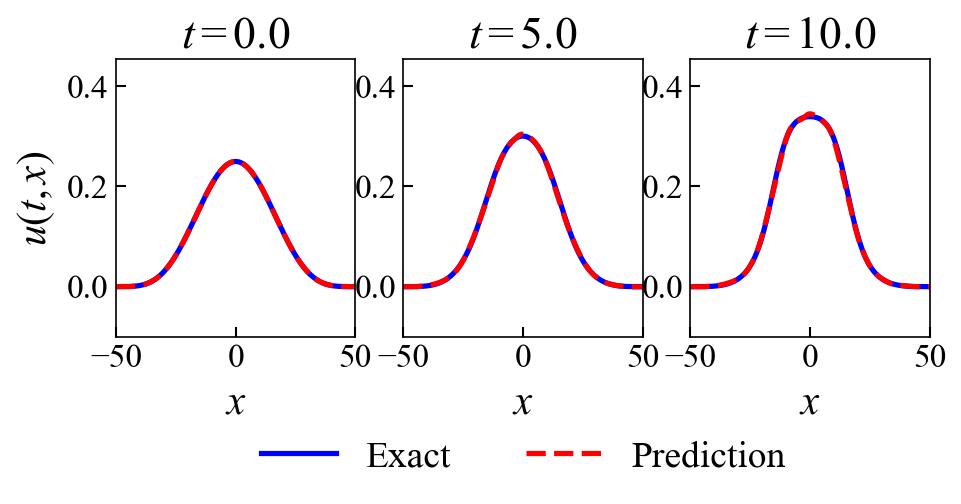}
    \end{minipage}
    \caption{Best solution of the  \textit{Gray--Scott}, obtained
with the configuration (PBC, IC2, SA, mini-batch). The resulting relative $L^2$ errors of u-species and v-species are $2.05\times10^{-4}$ and $8.46\times10^{-4}$ respectively. Full $(x,t)$  network predictions are in Fig.~\ref{fig:all}.}
    \label{fig:gs}
\end{figure}

Reaction and diffusion of chemical species can produce a variety of patterns \cite{de1991turing,nishiura1999skeleton}, reminiscent of those often seen in nature. The Gray--Scott type system is one of classical mathematical models for chemical reactions \cite{gray1983autocatalytic,gray1984autocatalytic,liang2022reversible}. The general irreversible Gray--Scott equations describe such reactions: \begin{equation} \label{eq:gs_reaction}
    \begin{aligned}
         U + 2V &\longrightarrow 3V, \\
         V &\longrightarrow P.
    \end{aligned}
\end{equation}
This system is defined by two equations that describe the dynamics of two reacting substances:
\begin{equation} \label{eq:GS}
    \begin{aligned}
u_t &= \epsilon_1 u_{xx} + b(1 - u) - uv^2, \\
v_t& = \epsilon_2  v_{xx} -(b+k)v + uv^2, \quad (t, x) \in (0, T]\times(-L, L), \\
u(0, x)& = u_0(x), \quad v(0, x) = v_0(x), \quad  x \in [-L, L], \\
u(t, -L)& = u(t, L), \quad v(t, -L) = v(t, L), \quad  t \in [0, T],
    \end{aligned}
\end{equation}
where $T=20,L=50$, $\epsilon_1=1,\epsilon_2=0.01$ are diffusion rates, $b=0.02$ is the ``feeding rate" that adds $U$, $k=0.0562$ is the  ``killing rate" that removes $V$. We set our initial conditions as:\begin{equation*}
    \begin{aligned}
        u_0(x) = 1- \frac{\sin(\pi (x - 50)/100)^4}{2}, \quad v_0(x) = \frac{\sin(\pi (x - 50)/100)^4}{4}.
    \end{aligned}
\end{equation*}
In this ablation study, we not only evaluate the performance of self-adaptive neural networks but also test our scheme with mini-batching on the Gray--Scott equations, which model the densities of two interacting species. For mini-batching, we use 5 batches, each containing 3,280 samples. As shown in Table \ref{tab:gs}, the results confirm that each proposed component improves the overall model performance. Omitting any one of these components leads to higher error rates. Notably, the errors for the $u-$equation are typically smaller than those for the  $v-$equation, which aligns with our expectations. This is because the $v-$equation is sharper and more challenging to train compared to the $u-$equation. Additionally, experiments using initial condition \textit{IC2} consistently outperform those using \textit{IC1}. Specifically, the full scheme with \textit{IC2} achieves the best results, with a relative  $L^2$ error  of $2.05\times 10^{-4}$ for $u-$equation and $8.46\times 10^{-4}$ for $v-$equation. These results suggest that \textit{IC2} is more effective for larger time-scale PDEs, as \textit{IC1} causes the initial condition to decay too quickly, reducing its influence during training. The predicted solutions from our best model are visualized in Fig.~\ref{fig:gs}.

% \begin{figure}[h]
% \centering
% \includegraphics[width=0.75\textwidth]{original_doc/fig_GS_u_BL+PBC+IC2+SA+minibatch.png}
% \caption{Best solution of the \textit{Gray-Scott  (u species)}. The resulting relative $L^2$ error is $2.05\times10^{-4}$. Full $(x,t)$  network predictions are in Fig.~\ref{fig:all}.}
% \label{fig:gs_u}
% \end{figure}
% \begin{figure}[h]
% \centering
% \includegraphics[width=0.75\textwidth]{original_doc/fig_GS_v_BL+PBC+IC2+SA+minibatch.png}
% \caption{Best solution of the  \textit{Gray-Scott (v species)}. The resulting relative $L^2$ error is $8.46\times10^{-4}$. Full $(x,t)$  network predictions are in Fig.~\ref{fig:all}.}
% \label{fig:gs_v}
% \end{figure}
% Figure \ref{fig:gs_utx} and Figure \ref{fig:gs_vtx} show the comparison  of the best prediction against the reference solution. Note that the majority of this error is attributed to the sharpness of the PDEs.
% \begin{figure}[h]
%  \centering
% \includegraphics[width=0.6\textwidth]{original_doc/fig3_GS_u_BL+PBC+IC2+SA+minibatch.png}
% \caption{Exact solution of the \textit{Gray-Scott  (v specie)} with the corresponding best network prediction and the absolute error difference. }
% \label{fig:gs_utx}
% \end{figure}

% \begin{figure}[h]
% \centering
% \includegraphics[width=0.6\textwidth]{original_doc/fig3_GS_v_BL+PBC+IC2+SA+minibatch.png}
% \caption{Exact solution of the  \textit{Gray-Scott  (v specie)} with the corresponding best network prediction and the absolute error difference.  }
% \label{fig:gs_vtx}
% \end{figure}

\begin{table}[ht]
\setlength{\tabcolsep}{0.8mm}
%title of the table
\centering
\renewcommand{\arraystretch}{1.5}
\begin{tabular}{ccccccc}
\toprule
\multicolumn{5}{c}{\textbf{Ablation Settings}}&\multicolumn{2}{c}{\textbf{Performance} }\\ [0.5ex]
\cmidrule(lr){1-5} \cmidrule(lr){6-7}

\textbf{PBC} & \textbf{IC1} & \textbf{IC2} & \textbf{SA}& \textbf{mini-batch} &\textbf{Rel. $L^2$ error}: $u$& \textbf{Rel. $L^2$ error}: $v$\\[0.2ex]
\hline

 % Entering row contents Midnight&7&-3&
\Checkmark & \Checkmark & \XSolidBrush &\Checkmark &\Checkmark & $8.98\times 10^{-4}$ & $7.92\times 10^{-3}$\\[0.5ex]
\Checkmark & \XSolidBrush & \Checkmark &\Checkmark & \Checkmark & $2.05\times 10^{-4}$ & $8.46\times 10^{-4}$\\[0.5ex]
\Checkmark & \Checkmark & \XSolidBrush &\XSolidBrush &\Checkmark &  $4.93\times 10^{-3}$ & $3.03\times 10^{-2}$\\[0.5ex]
\Checkmark & \XSolidBrush & \Checkmark &\XSolidBrush & \Checkmark & $2.86\times 10^{-3}$ & $8.20\times 10^{-3}$\\[0.5ex]
\Checkmark & \Checkmark & \XSolidBrush &\Checkmark &\XSolidBrush  &  $9.01\times 10^{-3}$ & $4.27\times 10^{-3}$\\[0.5ex]
\Checkmark & \XSolidBrush & \Checkmark & \Checkmark & \XSolidBrush &$7.59\times 10^{-3}$ & $3.98\times 10^{-3}$\\[0.5ex]
\Checkmark& \Checkmark & \XSolidBrush &\XSolidBrush & \XSolidBrush &$5.11\times 10^{-3}$ & $2.58\times 10^{-2}$\\[0.5ex]
\Checkmark& \XSolidBrush  & \Checkmark &\XSolidBrush &\XSolidBrush & $1.23\times 10^{-3}$ & $9.74\times 10^{-3}$\\[0.5ex]
\Checkmark& \XSolidBrush  & \XSolidBrush &\XSolidBrush & \XSolidBrush &$8.09\times 10^{-1}$ & $9.99\times 10^{-1}$\\[0.5ex]
\XSolidBrush  &  \Checkmark  & \XSolidBrush  &\XSolidBrush & \XSolidBrush &$2.06\times 10^{-1}$ & $7.43\times 10^{-1}$\\[0.5ex]
\XSolidBrush  &  \XSolidBrush  & \Checkmark  &\XSolidBrush & \XSolidBrush &$5.37\times 10^{-1}$ & $2.58\times 10^{-1}$\\[0.5ex]
\XSolidBrush  &  \XSolidBrush  & \XSolidBrush  &\Checkmark & \XSolidBrush &$1.15\times 10^{-1}$ & $2.62\times 10^{-1}$\\[0.5ex]
\XSolidBrush  &  \XSolidBrush  & \XSolidBrush  &\XSolidBrush & \Checkmark &$9.53\times 10^{-2}$ & $8.65\times 10^{-2}$\\[0.5ex]
\XSolidBrush  &  \XSolidBrush  & \XSolidBrush  &\XSolidBrush & \XSolidBrush &$9.21\times 10^{-2}$ & $2.98\times 10^{-1}$\\[0.5ex]
\bottomrule
\end{tabular}
\vspace{0.4cm}
\caption{\textit{Gray--Scott Equation:} Relative $L^2$ errors of $u$ and $v$ equations for an ablation study illustrating the impact of disabling individual components of the proposed technique and training pipeline.}
\label{tab:gs}
\end{table}

\subsection{Nonlinear Schr\"odinger equation}
In theoretical physics, nonlinear Schr\"odinger equation is a nonlinear PDE, applicable to classical and quantum mechanics \cite{fibich2015nonlinear, kato1987nonlinear,kevrekidis2001discrete}.  The dimensionless equation of the classical field is \begin{equation}
\begin{aligned}
     u_t &= i u_{xx} + i|u|^2u, \quad (t,x) \in [0, 2]\times[-\pi, \pi],  \\
     u(0, x) & = u_0(x), \quad  x \in [-\pi, \pi], \\
u(t, -\pi) &= u(t, \pi), \quad  t \in [0, 2],
\end{aligned}
\end{equation}
where we let $$u_0(x) =\frac{2}{2 - \sqrt{2} \cos(x)} - 1.$$

\begin{table}[h]
\setlength{\tabcolsep}{0.8mm}
%title of the table
\centering
\resizebox{\textwidth}{!}{
\renewcommand{\arraystretch}{1.6}
\begin{tabular}{ccccccc}
\toprule
\multicolumn{5}{c}{\textbf{Ablation Settings}}&\multicolumn{2}{c}{\textbf{Performance} }\\ [0.5ex]
\cmidrule(lr){1-5} \cmidrule(lr){6-7}

\textbf{PBC} & \textbf{IC1} & \textbf{IC2} & \textbf{SA}& \textbf{mini-batch} &\textbf{Rel. $L^2$ error}: \textit{real part}&\textbf{Rel. $L^2$ error}: \textit{imag part}\\[0.2ex]
\hline

 % Entering row contents Midnight&7&-3&
\Checkmark & \Checkmark & \XSolidBrush &\Checkmark &\Checkmark & $1.80\times 10^{-4}$ & $2.58\times 10^{-4}$\\[0.5ex]
\Checkmark & \XSolidBrush & \Checkmark &\Checkmark & \Checkmark & $1.84\times 10^{-4}$ & $3.27\times 10^{-4}$\\[0.5ex]
\Checkmark & \Checkmark & \XSolidBrush &\XSolidBrush &\Checkmark &  $7.34\times 10^{-3}$ & $1.23\times 10^{-2}$\\[0.5ex]
\Checkmark & \XSolidBrush & \Checkmark &\XSolidBrush & \Checkmark & $8.32\times 10^{-3}$ & $1.39\times 10^{-2}$\\[0.5ex]
\Checkmark & \Checkmark & \XSolidBrush &\Checkmark &\XSolidBrush  &  $1.10\times 10^{-3}$ & $2.10\times 10^{-3}$\\[0.5ex]
\Checkmark & \XSolidBrush & \Checkmark & \Checkmark & \XSolidBrush &$1.70\times 10^{-3}$ & $3.27\times 10^{-3}$\\[0.5ex]
\Checkmark& \Checkmark & \XSolidBrush &\XSolidBrush & \XSolidBrush &$1.00\times 10^{-1}$ & $1.83\times 10^{-1}$\\[0.5ex]
\Checkmark& \XSolidBrush  & \Checkmark &\XSolidBrush &\XSolidBrush & $9.56\times 10^{-2}$ & $1.56\times 10^{-1}$\\[0.5ex]
\Checkmark& \XSolidBrush  & \XSolidBrush &\XSolidBrush & \XSolidBrush &$8.09\times 10^{-1}$ & $9.99\times 10^{-1}$\\[0.5ex]
\XSolidBrush  &  \Checkmark  & \XSolidBrush  &\XSolidBrush & \XSolidBrush &$9.99\times 10^{-1}$ & $8.78\times 10^{-1}$\\[0.5ex]
\XSolidBrush  &  \XSolidBrush  & \Checkmark  &\XSolidBrush & \XSolidBrush &$5.04\times 10^{-2}$ & $1.21\times 10^{-1}$\\[0.5ex]
\XSolidBrush  &  \XSolidBrush  & \XSolidBrush  &\Checkmark & \XSolidBrush &$5.44\times 10^{-2}$ & $7.26\times 10^{-2}$\\[0.5ex]
\XSolidBrush  &  \XSolidBrush  & \XSolidBrush  &\XSolidBrush & \Checkmark &$6.60\times 10^{-2}$ & $1.05\times 10^{-1}$\\[0.5ex]
\XSolidBrush  &  \XSolidBrush  & \XSolidBrush  &\XSolidBrush & \XSolidBrush &$3.63\times 10^{-1}$ & $5.06\times 10^{-1}$\\[0.5ex]

\bottomrule
\end{tabular}}
\vspace{0.4cm}
\caption{\textit{Nonlinear Schr\"odinger Equation:} Relative $L^2$ errors of real and imaginary equations for an ablation study illustrating the impact of disabling individual components of the proposed technique and training pipeline.}
\label{tab:nls}
\end{table}

Using the same parameter settings and training steps with prior experiments, we obtain the following results:

\begin{figure}[ht]
    \centering
    % First figure
    \begin{minipage}{0.48\textwidth}
        \centering
        \includegraphics[width=\linewidth]{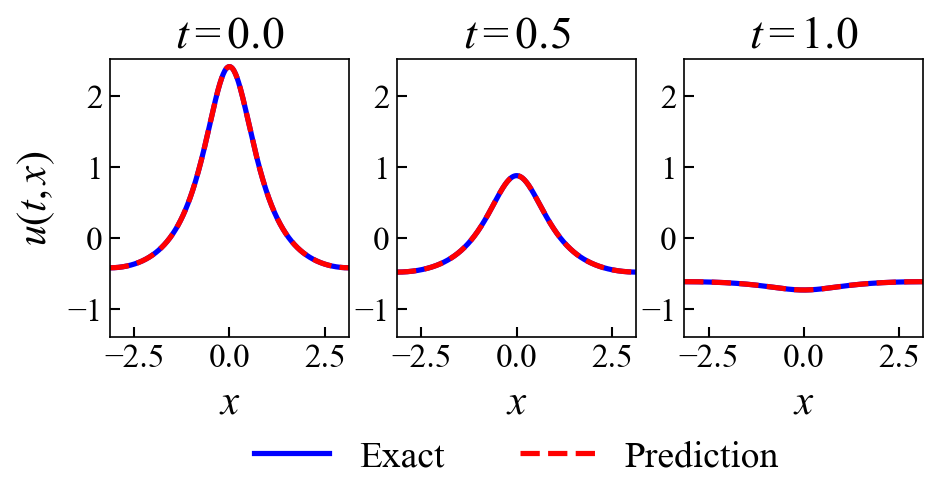}
    \end{minipage}
    \hfill
    % Second figure
    \begin{minipage}{0.48\textwidth}
        \centering
        \includegraphics[width=\linewidth]{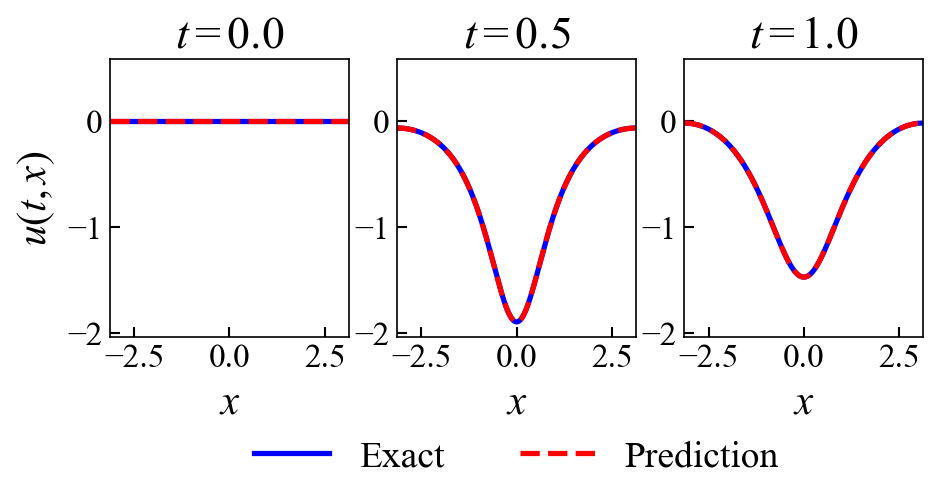}
    \end{minipage}
    \caption{Best solution of the \textit{Nonlinear Schr\"odinger Equation},
    obtained with the configuration (PBC, IC1, SA,
    mini-batching). The resulting relative $L^2$ errors of real-part and
    imaginary-part solutions are $1.80\times10^{-4}$ and
   $2.58\times10^{-4}$ respectively. Full $(x,t)$ network
    predictions are in Fig.~\ref{fig:all}.}
    \label{fig:nls}
\end{figure}

Table \ref{tab:nls} presents the results of all experiments. The full scheme achieves the best performance, with a relative $L^2$ error of real part of $1.80\times 10^{-4}$ and the relative $L^2$ error of imaginary part of $2.58\times 10^{-4}$. There is no significant difference in the errors between methods using \textit{IC1} and those using \textit{IC2}. Moreover, without the use of self-adaptive neural networks or mini-batching, accurately obtaining the solution becomes challenging, highlighting the importance of both methods in ensuring more stable training.

\subsection{Linear advection}\label{sec:la}
The linear advection equation models the transport of a quantity by a constant velocity field and serves as a fundamental benchmark for testing numerical methods due to its simplicity and rich structure \cite{leveque2002finite, toro2013riemann}. It describes the propagation of a scalar profile without deformation and is given by:
\begin{equation}
\begin{aligned}
u_t + c u_x &= 0, \quad\quad\quad (t, x) \in (0, T] \times (0, 2\pi), \\
u(0, x) &= \sin (x), \quad x \in [0, 2\pi], \\
u(t, 0) &= u(t, 2\pi), \quad t \in [0, T].
\end{aligned}\label{eq:la}
\end{equation}
Unlike dissipative PDEs such as the Allen–Cahn or Cahn–Hilliard equations, the linear advection equation conserves the shape of the initial profile as it translates. This problem becomes ``stiffer” as the advection velocity increases.

% \begin{theorem}\label{thm:advection}
% Let $u(t,x)$ denote the solution to the linear advection equation
% \eqref{eq:la}. With the transformed neural approximation
% {$\tilde u = \psi + \phi\,u_{nn}$ defined in
% \eqref{eq:trans}, where $\psi$ and $\phi$ satisfy the conditions in
% \eqref{eq:HC_condition}}, $u=\tilde u$ satisfies \eqref{eq:la} if and
% only if $u_{nn}$ satisfies the forced linear transport
% equation
% \[
% \phi\,\bigl(\partial_t u_{nn} + c\,\partial_x u_{nn}\bigr)
% + \bigl(\phi_t + c\,\phi_x\bigr)\,u_{nn}
% = -\bigl(\psi_t + c\,\psi_x\bigr).
% \]
% This transformed problem is well-posed and stable for PINN training
% provided the forcing $-(\psi_t+c\psi_x)\in C^1([0,T]\times\Omega)$ and is
% uniformly bounded or decaying in time.
% \end{theorem}
\begin{table}[h!]
\centering
\renewcommand{\arraystretch}{1.5}
\begin{tabular}{@{}llcc@{}}
\toprule
\textbf{Solver} & $c$ &  \textbf{Rel. $L^1$ error} &  \textbf{Rel. $L^2$ error}  \\
\midrule
 & 30 & $0.0030$ & $0.0041$  \\
HC-PINN    & 40 & $0.0141$& $0.0159$ \\
     & 50 & $0.0147$ & $0.0164$  \\
     & 60 & $0.0171$ &$0.0182$  \\
\midrule
 & 30 & $0.0478$ & $0.0579$  \\
 CINN\cite{braga2022characteristics}    & 40 & $0.0714$ & $0.0852$ \\
   & 50 & $0.4692$ & $0.5365$ \\
     & 60 & $0.7293$ & $0.7134$ \\
     \midrule
 & 30 & $0.0873$ & $0.1003$  \\
  Baseline    & 40 & $0.3918$ & $0.4395$ \\
PINN    & 50 & $0.7140$ & $0.7797$ \\
     & 60 & $0.8339$ & $0.8664$ \\
\bottomrule
\end{tabular}
\vspace{0.4cm}
\caption{Comparison of HC-PINN, CINN and baseline PINN across different advection speeds $c$  over 10 independent repetitions.} \label{tab:la}
\end{table}

\begin{theorem}[Transformation for the linear advection equation]
\label{thm:advection}
Consider the linear advection problem \eqref{eq:la}. Let
\[
    \tilde u(t,x)=\psi(t,x)+\phi(t,x)u_{nn}(t,x)
\]
be the transformed neural approximation defined in \eqref{eq:trans}, where
$\psi$ and $\phi$ satisfy the initial-condition requirements in
\eqref{eq:HC_condition}. Assume further that $\psi$, $\phi$, and $u_{nn}$
are periodic in $x$ with period $2\pi$. Then $\tilde u$ satisfies
\eqref{eq:la} if and only if $u_{nn}$ satisfies
\[
    \phi\bigl(\partial_t u_{nn}+c\,\partial_x u_{nn}\bigr)
    +
    \bigl(\phi_t+c\,\phi_x\bigr)u_{nn}
    =
    -\bigl(\psi_t+c\,\psi_x\bigr).
\]
This transformed problem is well-posed and stable for PINN training
provided the forcing $-(\psi_t+c\psi_x)\in C^1([0,T]\times\Omega)$ and is
uniformly bounded or decaying in time.
\end{theorem}

\begin{proof}
By the transformation \eqref{eq:trans},
\[
\begin{aligned}
    \partial_t\tilde u+c\,\partial_x\tilde u
    &=
    \partial_t(\psi+\phi u_{nn})
    +
    c\,\partial_x(\psi+\phi u_{nn}) \\
    &=
    \psi_t+c\,\psi_x
    +
    \phi\bigl(\partial_t u_{nn}+c\,\partial_x u_{nn}\bigr)
    +
    \bigl(\phi_t+c\,\phi_x\bigr)u_{nn}.
\end{aligned}
\]
Therefore,
\[
    \partial_t\tilde u+c\,\partial_x\tilde u=0
\]
if and only if
\[
    \phi\bigl(\partial_t u_{nn}+c\,\partial_x u_{nn}\bigr)
    +
    \bigl(\phi_t+c\,\phi_x\bigr)u_{nn}
    =
    -\bigl(\psi_t+c\,\psi_x\bigr).
\]
Moreover, the assumptions in \eqref{eq:HC_condition} imply
\[
    \tilde u(0,x)
    =
    \psi(0,x)+\phi(0,x)u_{nn}(0,x)
    =
    u_0(x),
\]
and the periodicity of $\psi$, $\phi$, and $u_{nn}$ implies
\[
    \tilde u(t,0)=\tilde u(t,2\pi).
\]
Hence $\tilde u$ satisfies the complete initial-boundary value problem
\eqref{eq:la} if and only if $u_{nn}$ satisfies the transformed equation
above.
\end{proof}

For our implementation, we take \[ \psi(t,x)=u_0(x)e^{-\alpha t}, \qquad \phi(t)=\frac{t}{T}, \qquad \alpha>0. \]
Since $u_0(x)=\sin x$ is periodic, both $\psi$ and $\tilde u$ are periodic in $x$ under the periodic neural network construction for $u_{nn}$. Moreover, \[ \tilde u(0,x)=u_0(x), \qquad \phi_t(0,x)=\frac{1}{T}\neq 0, \] so the transformation is non-degenerate at $t=0$ in the sense of Theorem~\ref{thm:wellposed}. Since $\phi_t=1/T$ and $\phi_x=0$, the transformed equation becomes \[ \frac{t}{T} \bigl( \partial_tu_{nn}+c\,\partial_xu_{nn} \bigr) + \frac{1}{T}u_{nn} = -\bigl(\psi_t+c\,\psi_x\bigr). \] With $u_0(x)=\sin x$, the right-hand side is \[ -\bigl(\psi_t+c\,\psi_x\bigr) = \bigl(\alpha\sin x-c\cos x\bigr)e^{-\alpha t}, \] which is smooth and decays exponentially in time,  guaranteeing both the stability and the regularity
of the learned component $u_{nn}(t,x)$.
% {For our implementation we take $\psi(t,x)=u_0(x)\,
% e^{-\alpha t}$ and $\phi(t)=t/T$ with $\alpha>0$, so that
% $\tilde u(0,x)=u_0(x)$ holds exactly and $\phi_t(0,x)=1/T\neq 0$,
% ensuring the transformation is non-degenerate at $t=0$ in the sense of
% Theorem~\ref{thm:wellposed}. Substituting $\phi=t/T$ (so $\phi_t=1/T$,
% $\phi_x=0$) yields the equation satisfied by $u_{nn}$,}
% \[
% {(t/T)\,\bigl(\partial_t u_{nn}
% + c\,\partial_x u_{nn}\bigr) + (1/T)\,u_{nn}
% = -\bigl(\psi_t + c\,\psi_x\bigr).}
% \]
% {With $u_0(x)=\sin x$ the right-hand side evaluates to
% $(\alpha\sin x - c\cos x)\,e^{-\alpha t}$, which is smooth and decays
% exponentially in time, guaranteeing both the stability and the regularity
% of the learned component $u_{nn}(t,x)$.}
Here, we compare the results
produced by our method with those from a baseline PINN under the same
hyperparameter settings. In our experiments, both models were able to
converge at low velocities; however, we observed that higher advection
velocities introduced stiffness. Our method exhibited greater robustness
and was less affected by stiffness compared to the baseline PINN.

\begin{figure}[h!]
 \centering
\includegraphics[width=0.85\textwidth]{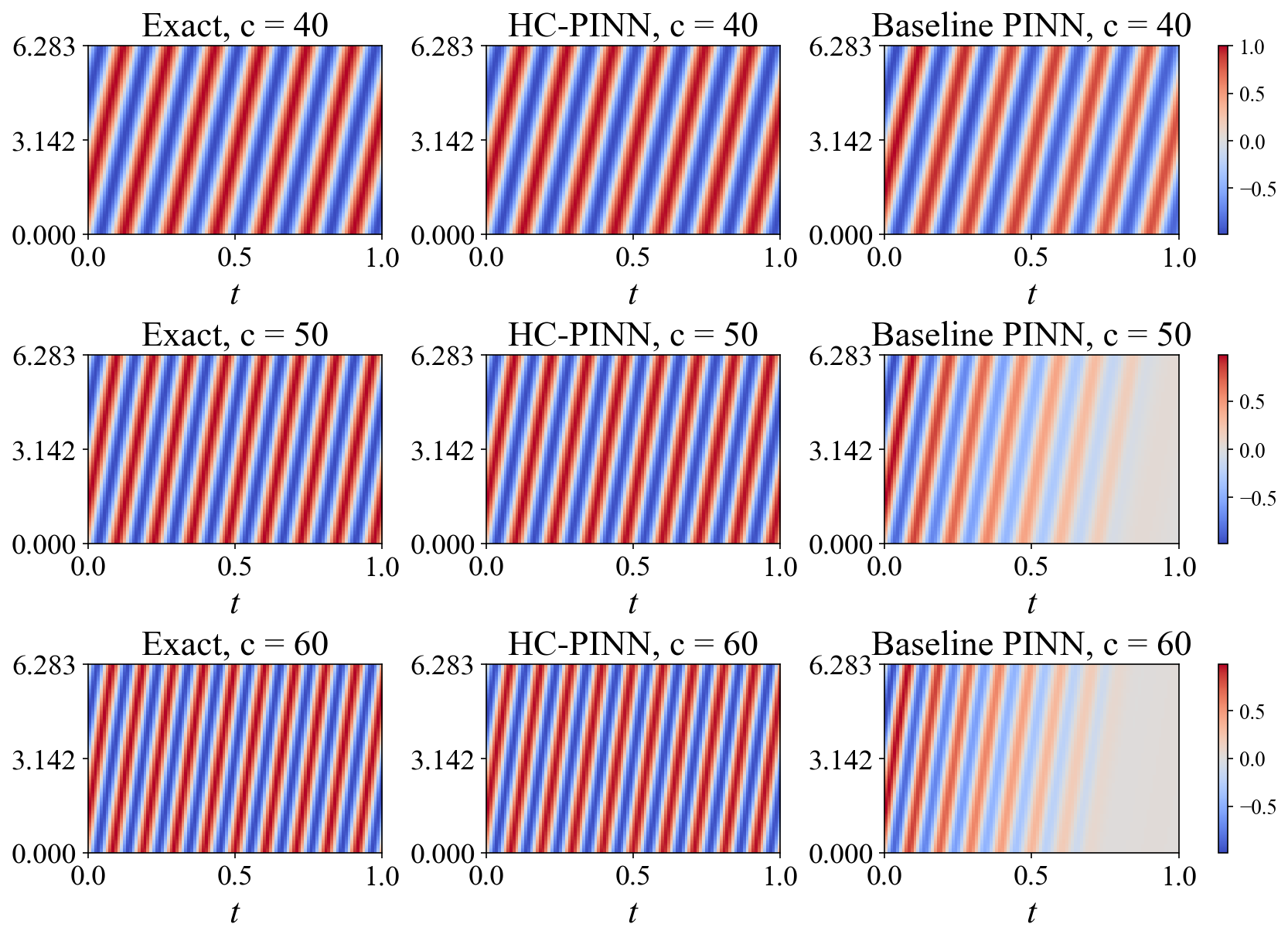}
% \caption{Periodic advection results after 50,000 iterations.}
\caption{{Periodic advection results after 50,000 iterations:
exact solution, HC-PINN, and baseline PINN for $c=40,50,60$ (cf.
Table~\ref{tab:la}).}}
\label{fig:la}
\end{figure}

\begin{table}[h]
%title of the table
\centering
\renewcommand{\arraystretch}{1.3}
\begin{tabular}{cccccc}
\toprule
\multicolumn{4}{c}{\textbf{Ablation Settings}}&\multicolumn{2}{c}{\textbf{Performance} }\\ [0.5ex]
\cmidrule(lr){1-4} \cmidrule(lr){5-6}

\textbf{PBC} & \textbf{IC} & \textbf{SA}& \textbf{Rel. $L^2$ error}& \textbf{Rel. $L^1$ error}\\[0.2ex]
\hline
\Checkmark & \Checkmark &\Checkmark & $1.82\times 10^{-2}$ & $1.71\times 10^{-2}$\\[0.5ex]
\Checkmark  & \Checkmark & \XSolidBrush& $9.75\times 10^{-2}$ & $8.65\times 10^{-2}$ \\[0.5ex]
\Checkmark &  \XSolidBrush &\Checkmark & $7.50\times 10^{-1}$ & $7.87\times 10^{-1}$\\[0.5ex]
\XSolidBrush & \Checkmark & \Checkmark & $3.87\times 10^{-2}$ & $5.47\times 10^{-2}$\\[0.5ex]
 \Checkmark & \XSolidBrush  &\XSolidBrush & $7.44\times 10^{-1}$ & $8.01\times 10^{-1}$\\[0.5ex]
\XSolidBrush &   \Checkmark &\XSolidBrush & $1.07\times 10^{-1}$ & $1.19\times 10^{-1}$\\[0.5ex]
\XSolidBrush  &  \XSolidBrush  &\Checkmark & $7.65\times 10^{-1}$ & $7.71\times 10^{-1}$\\[0.5ex]
\XSolidBrush  &  \XSolidBrush  & \XSolidBrush & $8.66\times 10^{-1}$ & $8.34\times 10^{-1}$\\[0.5ex]
\bottomrule
\end{tabular}
\vspace{0.4cm}
\caption{\textit{Linear Advection Equation (c=60):} Relative $L^2$ and $L^1$ errors for an ablation study illustrating the impact of disabling individual components of the proposed technique and training pipeline.} \label{tab:la-abl}
\end{table}

% For our choice, we have $\psi(t,x) = u_0(x)e^{-\alpha t}$ and $\phi(t,x) = t$ with $\alpha>0$, ensuring $\eta(0,x) = 0$. Then, the forcing term satisfies $\eta_t + c\eta_x = (\alpha \sin{x}-c\cos{x})e^{-\alpha t}$. Since this forcing is smooth and decays exponentially in time, it guarantees both the stability and regularity of the neural network $\eta(t,x)$.  Here, we compare the results produced by our method with those from a baseline PINN under the same hyperparameter settings. In our experiments,  both models were able to converge at low velocities; however, we observed that higher advection velocities introduced stiffness.  Our method exhibited greater robustness and was less affected by stiffness compared to the baseline PINN.

The implementation details are the same as in Section~\ref{section:B}, except that we perform $10$ independent replications of the experiment to account for the randomness in neural network weight initialization, rather than using the same initial weights. Table~\ref{tab:la} shows the averaged relative $L^1$/$L^2$-error over $10$ independent repetitions.
{As the advection velocity increases, the PDE becomes
stiffer and harder for the baseline PINN to train. Our method, however,
consistently produces better results, as clearly reflected in
Table~\ref{tab:la}.}
% As the advection velocity increases and the PDE becomes stiffer, making it harder for baseline PINN to train.
% However, our method consistently produces better results. This trend is clearly reflected in Table~\ref{tab:la}.
Even at low velocities, our method demonstrates significantly better accuracy compared to the baseline PINN. As $c$ increases, both the relative $L^1$ and $L^2$ errors grow, however, the errors associated with the PINN rise significantly faster than those of our method. The most striking example occurs at $c=50$, where the {baseline PINN} fails to approximate the solution, while our method still produces stable and good results. This behavior is evident in Table~\ref{tab:la} and Figure~\ref{fig:la}. Table~\ref{tab:la-abl} shows the ablation study of linear advection equation.

\subsection{ Comparison to other {state-of-the-art } PINNs }
Figure~\ref{fig:all} provides a summary of our ablation study over different benchmarks; it shows the comparison between the reference exact solutions and our best HC-PINN predictions, together with the corresponding absolute errors; the benchmarks include Allen--Cahn, Cahn--Hilliard, Kuramoto--Sivashinsky, Gray--Scott, and nonlinear Schr\"odinger. Complementing these visual comparisons, Table~\ref{tab:comparison} summarizes the comparison study between HC-PINN and two other PINN variants in terms of relative $L^2$ errors for Allen--Cahn, Cahn--Hilliard and Kuramoto--Sivashinsky equations under a strictly controlled setting: {all methods are trained using the same
optimization budget: \(50{,}000\) Adam iterations followed by \(5{,}000\)
L-BFGS-B iterations,} same data and parameter settings to ensure fairness, unless otherwise specified.  {We emphasize that this controlled
setting is intended to isolate the effect of the training scheme under an
equal compute budget, and we note that the shorter budget may be
unfavorable to the baseline methods; for transparency, the originally
reported $300k$-iteration values of the baselines are included in
Table~\ref{tab:comparison} as a reference. Under this equal-budget protocol, HC-PINN attains substantially lower errors across all three benchmarks (while at full 300k budget Enhanced RBA matches or surpasses HC-PINN on AC Case $I$, our method remains substantially more accurate at the same compute budget).} Overall, these extensive experiments highlight the effectiveness and adaptability of our PINN scheme in solving a wide variety of stiff time-dependent PDEs. The results provide not only strong empirical validation but also a reliable and reproducible benchmark for future developments in this area.

\begin{table}[H]
\centering
\resizebox{\textwidth}{!}{
\renewcommand{\arraystretch}{1.5}
\begin{tabular}{ c c c c c c }
\toprule
  & \multicolumn{2}{c}{{Original ($300k$)}}
  & \multicolumn{3}{c}{{Equal budget ($50k$)}} \\
\cmidrule(lr){2-3}\cmidrule(lr){4-6}
  & {Causal} & {RBA}
  & Causal (MLP) &  RBA & HC-PINNs \\ \hline
 AC (Case $I$)  & {$1.43 \times 10^{-3}$} & {$5.7 \times 10^{-5}$} & $6.95\times 10^{-2}$ & $2.62\times 10^{-3}$ & $8.29\times 10^{-4}$ \\
 AC (Case $II$) & {$-$} & {$-$} & $1.78\times 10^{-2}$ & $1.08\times 10^{-3}$ & {$3.53\times 10^{-4}$} \\
 CH           & {$-$} & {$-$} & $3.49\times 10^{-1}$ & $9.83\times 10^{-2}$ & {$9.75\times 10^{-4}$} \\
 KS           & {$2.26 \times 10^{-2}$} & {$-$} & $2.72\times 10^{-1}$ & $3.64\times 10^{-2}$ & $5.37\times 10^{-4}$ \\
\bottomrule
\end{tabular}}
\caption{Relative $L^2$ errors for the Allen--Cahn (AC), Cahn--Hilliard (CH)
and Kuramoto--Sivashinsky (KS) equations. {The right block
reports results under a strictly controlled equal-budget protocol in which
all methods are trained for $50k$ iterations  followed by
\(5{,}000\) L-BFGS-B iterations with the same data and
parameter settings. The left block reports the originally published $300k$-iteration values of Causal training~\cite{wang2022respecting} and Enhanced RBA~\cite{anagnostopoulos2024residual}, included as a transparent reference. Entries marked $-$ indicate that the corresponding benchmark or case was not reported in the original paper at $300k$ iterations. All HC-PINN entries are relative $L^2$ errors of the best-performing configuration: AC~(Case~$I$) and AC~(Case~$II$) use (PBC,~IC1,~SA), CH uses (PBC,~IC2,~SA), and KS uses (PBC,~IC2,~SA), consistent with the corresponding ablation tables and best-solution figures.}}
\label{tab:comparison}
\end{table}

\begin{figure}[h]
    \centering
    \makebox[0.27\textwidth][c]{\small\textbf{Exact $u(x,y)$}}%
    \hfill
    \makebox[0.27\textwidth][c]{\small\textbf{Predicted $u(x,y)$}}%
    \hfill
    \makebox[0.25\textwidth][c]{\small\textbf{Error}}

    \vspace{0.3em}

    % ====== Row 1 ======
     \begin{minipage}[c]{0.03\textwidth} % left label
        \centering \textbf{(a)}
    \end{minipage}
    \begin{minipage}[c]{0.27\textwidth}
        \centering
        \includegraphics[width=\linewidth]{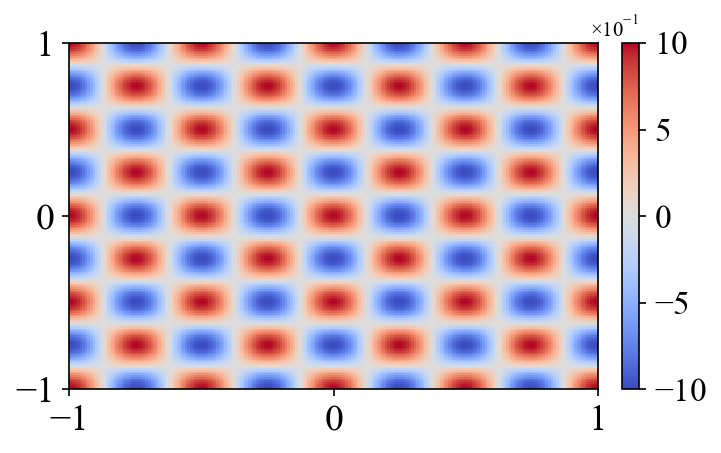}
        %\captionof*{figure}{Exact $t=0$}

    \end{minipage}
    \hfill
    \begin{minipage}[c]{0.27\textwidth}
        \centering
        \includegraphics[width=\linewidth]{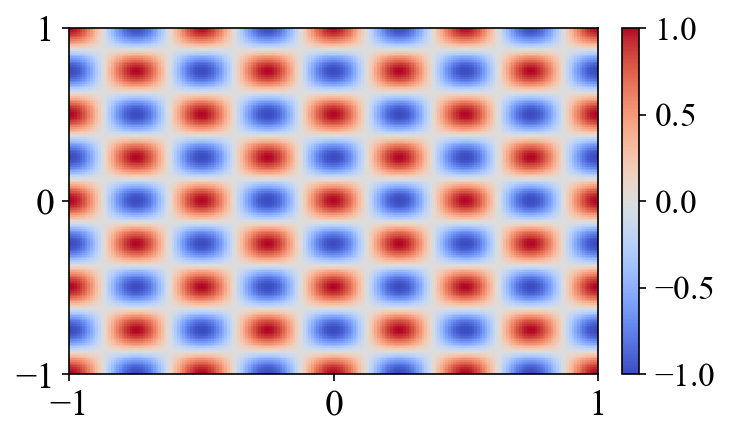}
       % \captionof*{figure}{Predicted $t=0$}
    \end{minipage}
    \hfill
    \begin{minipage}[c]{0.25\textwidth}
        \centering
        \includegraphics[width=\linewidth]{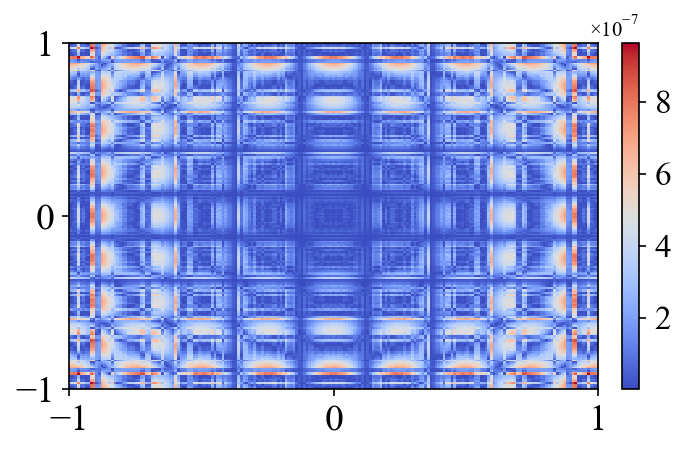}
       % \captionof*{figure}{Error $t=0$}
    \end{minipage}

    \vspace{0.5em}

    % ====== Row 2 ======
     \begin{minipage}[c]{0.03\textwidth} % left label
       \centering \textbf{(b)}
    \end{minipage}
    \begin{minipage}[c]{0.27\textwidth}
        \centering
        \includegraphics[width=\linewidth]{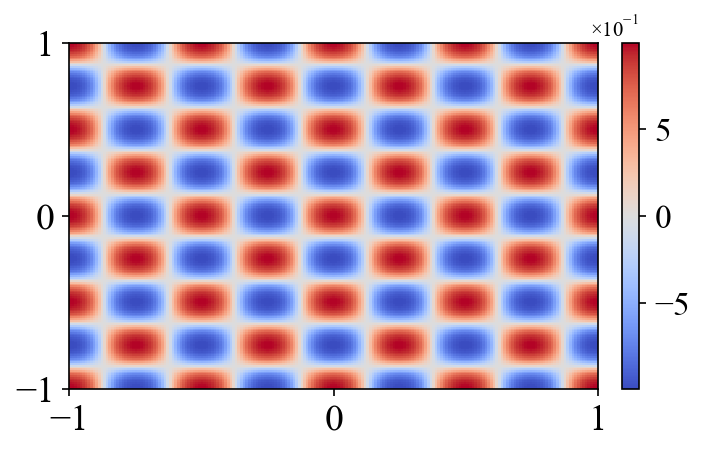}
        %\captionof*{figure}{Exact $t=0.5$}
    \end{minipage}
    \hfill
    \begin{minipage}[c]{0.27\textwidth}
        \centering
        \includegraphics[width=\linewidth]{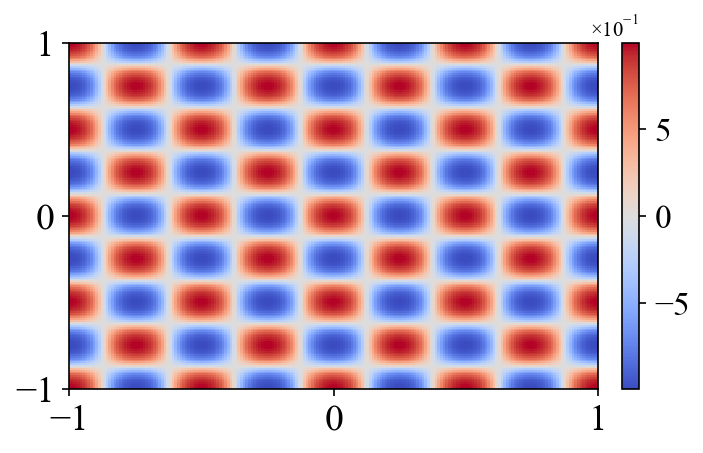}
        %\captionof*{figure}{Predicted $t=0.5$}
    \end{minipage}
    \hfill
    \begin{minipage}[c]{0.25\textwidth}
        \centering
        \includegraphics[width=\linewidth]{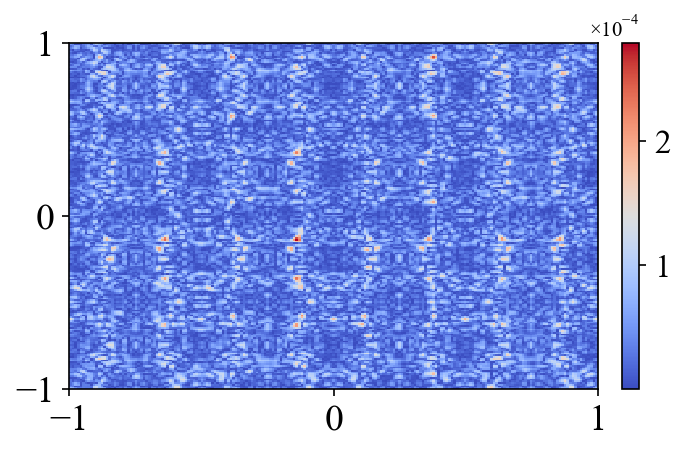}
       % \captionof*{figure}{Error $t=0.5$}
    \end{minipage}

    \vspace{0.5em}

    % ====== Row 3 ======
    \begin{minipage}[c]{0.03\textwidth} % left label
       \centering \textbf{(c)}
    \end{minipage}
    \begin{minipage}[c]{0.27\textwidth}
        \centering
        \includegraphics[width=\linewidth]{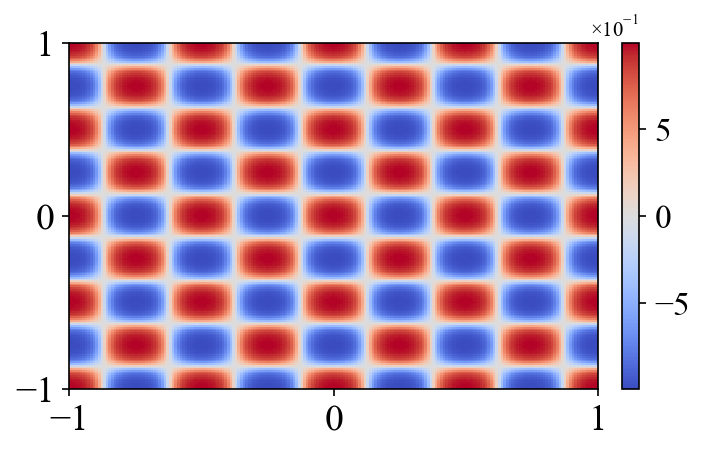}
        %\captionof*{figure}{Exact $t=1$}
    \end{minipage}
    \hfill
    \begin{minipage}[c]{0.27\textwidth}
        \centering
        \includegraphics[width=\linewidth]{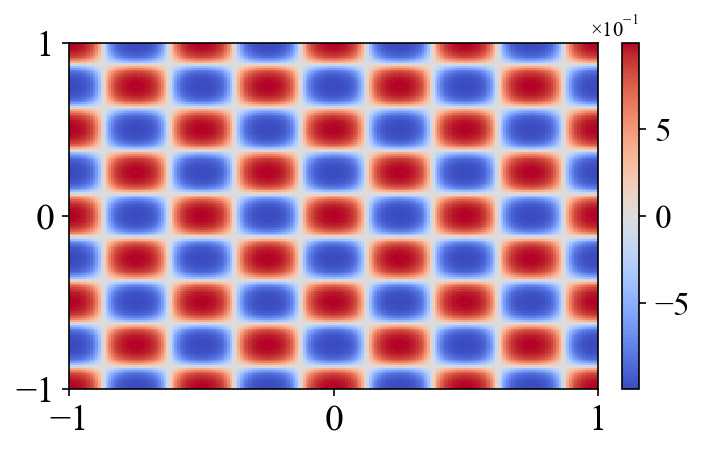}
    \end{minipage}
    \hfill
    \begin{minipage}[c]{0.25\textwidth}
        \centering
        \includegraphics[width=\linewidth]{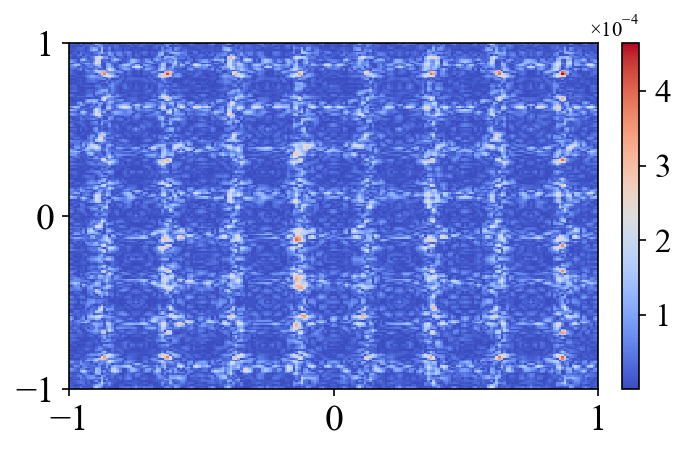}
    \end{minipage}

    \caption{Reference, HC-PINN predicted solution, and absolute error of the $2D$ Allen Cahn equation at different time
snapshots (a) $t = 0$, (b) $t = 0.5$, (c) $t = 1$.}\label{2d-allen}
\end{figure}

\subsection{Two-dimensional Allen--Cahn equation and sensitivity test}
The results of solving the $2D$ Allen–Cahn equation in Figure~\ref{2d-allen} were obtained after training the network for $50k$ iterations using the initial condition $u(0,x,y) = $\break$\sin(4\pi x)\cos(4\pi y)$, with $\gamma_1 = 0.0001$ and $\gamma_2 = 1$ in Eq.~\ref{eq:ac}, and periodic boundary conditions. The HC-PINN accurately captures the solution evolution, yielding exceptionally small relative $L^2$ errors: $8.04\times10^{-5}$ at $t=0.5$ and $1.33\times10^{-4}$ at $t=1$. This demonstrates the high accuracy and stability of our method for 2D Allen–Cahn dynamics. At last, we perform a sensitivity study on the number of Fourier modes $m$ for the Allen–Cahn benchmarks to assess how this key hyperparameter influences the accuracy of HC-PINNs. The results demonstrate that accuracy is highly dependent on spectral resolution: while small $m$ yields stable and low errors, performance deteriorates significantly once $m$ exceeds a critical threshold.

\begin{figure}[h!]
    \centering
    % First figure
    \begin{minipage}{0.45\textwidth}
        \centering
        \includegraphics[width=\linewidth]{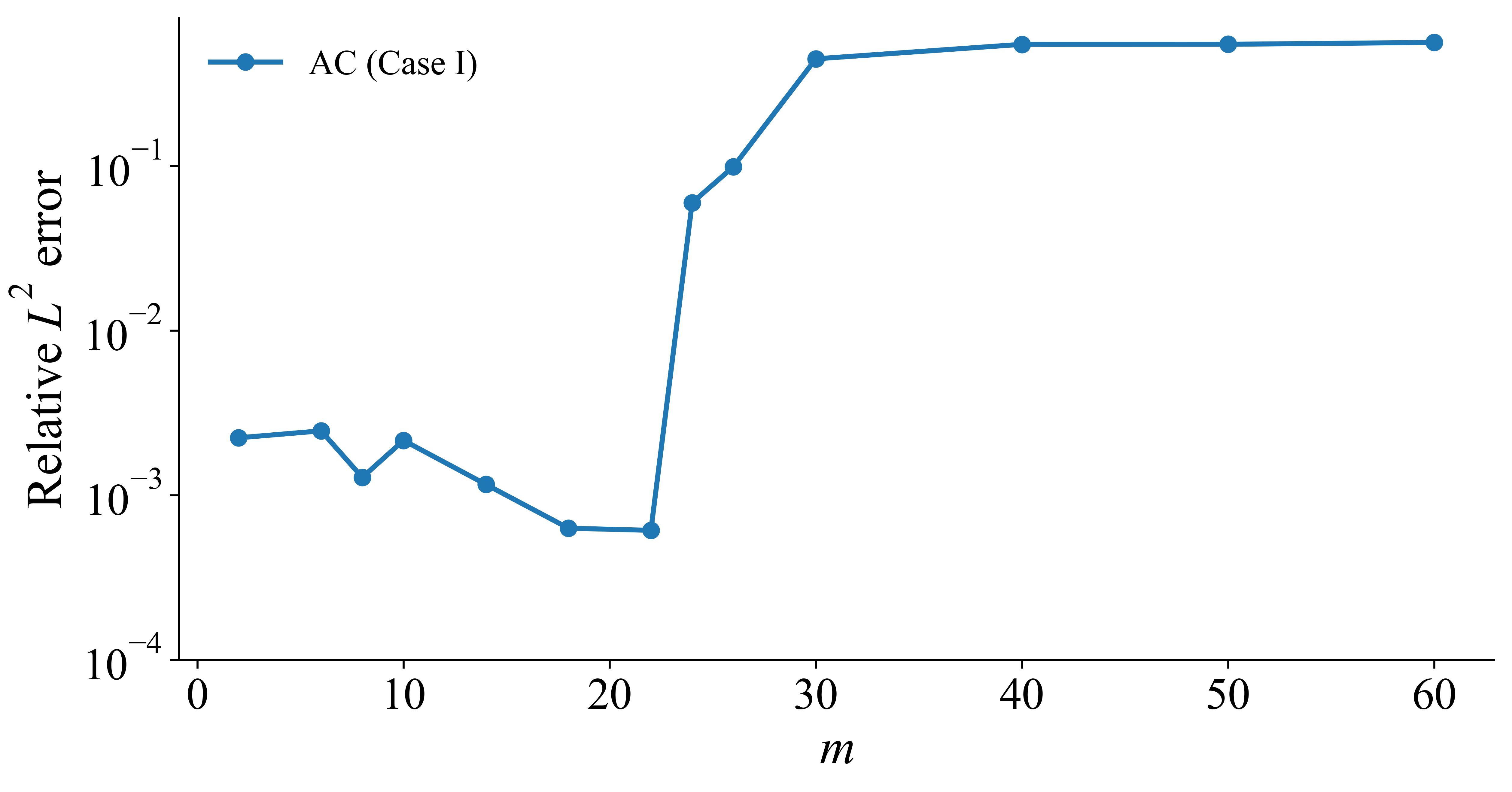}
    \end{minipage}
    \hfill
    % Second figure
    \begin{minipage}{0.45\textwidth}
        \centering
        \includegraphics[width=\linewidth]{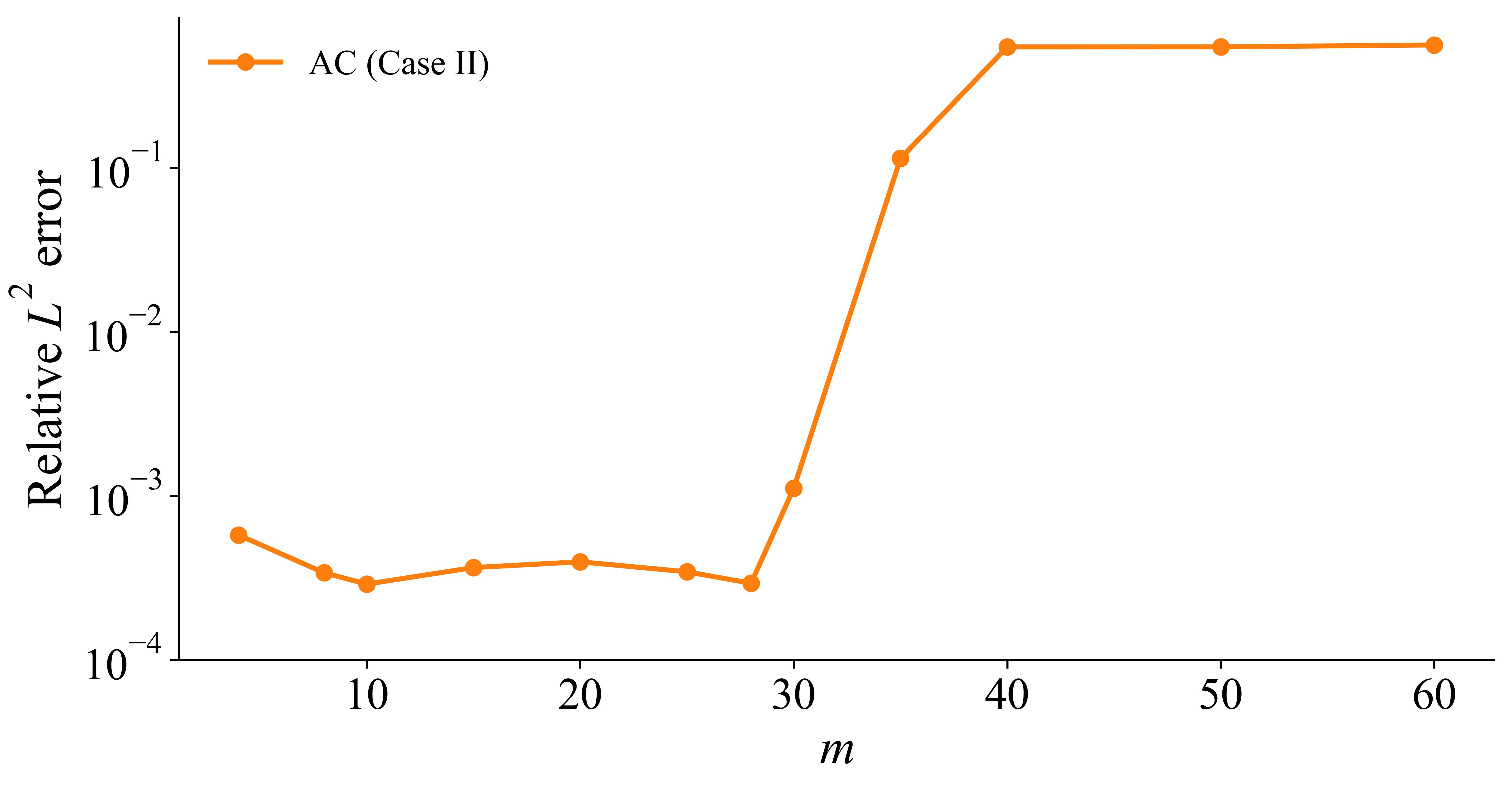}
    \end{minipage}

    \caption{Sensitivity of HC-PINNs to the number of Fourier modes
$m$ in the Allen–Cahn equations. The plot shows the relative
$L^2$ error as a function of $m$.}
    \label{fig:errorvsm}
\end{figure}

\subsection{{Robustness to non-periodic boundary conditions}}
\label{sec:dirichlet}

{All examples in Sections~\ref{sec:ac}--\ref{sec:la} use
periodic boundary conditions and rely on the Fourier-feature architecture
to enforce periodicity in $u_{nn}$. To assess whether the stabilizing
effect documented in those sections is intrinsic to the
hard-initial-condition mechanism rather than to the periodic
architecture, we test HC-PINN on the same Allen--Cahn problem of
Case~$I$ (Section~\ref{sec:ac}) but with non-homogeneous Dirichlet
boundary conditions:}
\begin{equation}\label{eq:ac_dirichlet}
{\begin{aligned}
u_t - \gamma_1 u_{xx} + \gamma_2(u^3 - u) & = 0, \quad (t,x)\in(0,1]\times(-1,1),\\
u(0,x) &= x^2\cos(\pi x), \quad x\in[-1,1],\\
u(t,-1) &= u(t,1) = -1, \quad t\in[0,1],
\end{aligned}}
\end{equation}
{with $\gamma_1 = 10^{-3}$ and $\gamma_2 = 5$ as in
Case~$I$. The boundary value $g = -1$ is dictated by the compatibility
condition $u_0(\pm 1) = -1$, so the same initial profile from the
periodic Case~$I$ is reused without modification.}

{The general hard-constraint construction~\eqref{eq:trans}
admits the boundary-aware lift}
\[
{\psi(t,x) = u_0(x)\,e^{-t} + g\,\bigl(1 - e^{-t}\bigr),
\qquad
\phi(t,x) = t\,(1 - x^2),}
\]
{which satisfies $\psi(0,x)=u_0(x)$, $\psi(t,\pm 1)=g$
for all $t$, $\phi(0,x)=0$, and $\phi(t,\pm 1)=0$. Hence
$\tilde u(0,x)=u_0(x)$ and $\tilde u(t,\pm 1)=g$ are enforced
exactly by construction, irrespective of $u_{nn}$. In particular,
neither an initial-condition loss nor a boundary-condition loss is
needed, and the training reduces to a single PDE-residual objective.}

{Crucially, in this experiment $u_{nn}(t,x)$ is a
standard fully connected MLP taking $(t,x)$ directly as input,
without any Fourier-feature layer. The architecture, optimizer
schedule, and collocation budget are identical to the periodic
experiments of Section~\ref{section:B}.}

\begin{figure}[h]
    \centering
    \includegraphics[width=0.8\textwidth]{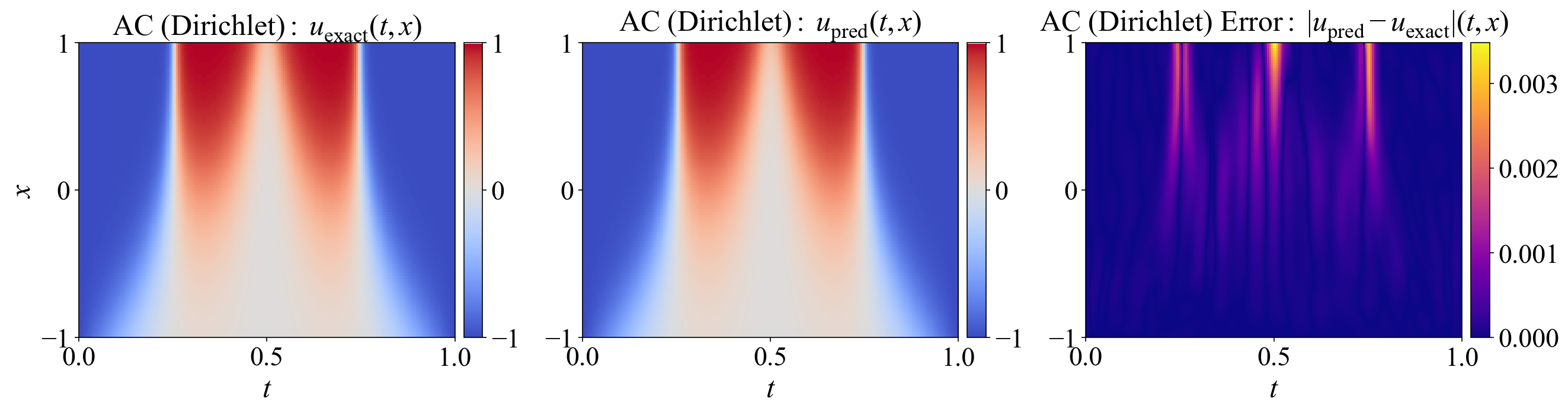}
    \caption{{Allen--Cahn with non-homogeneous Dirichlet
    BC ($g=-1$), $\gamma_1=10^{-3}$, $\gamma_2=5$. Left: reference
    solution by Crank--Nicolson IMEX FD. Middle: HC-PINN prediction with
    plain MLP (no Fourier features) and hard IC/BC enforcement. Right:
    pointwise absolute error. Relative $L^2$ error $5.05\times10^{-4}$,
    relative $L^1$ error $2.82\times10^{-4}$.}}
    \label{fig:ac_dirichlet}
\end{figure}

\begin{figure}[h]
    \centering
    \includegraphics[width=0.8\textwidth]{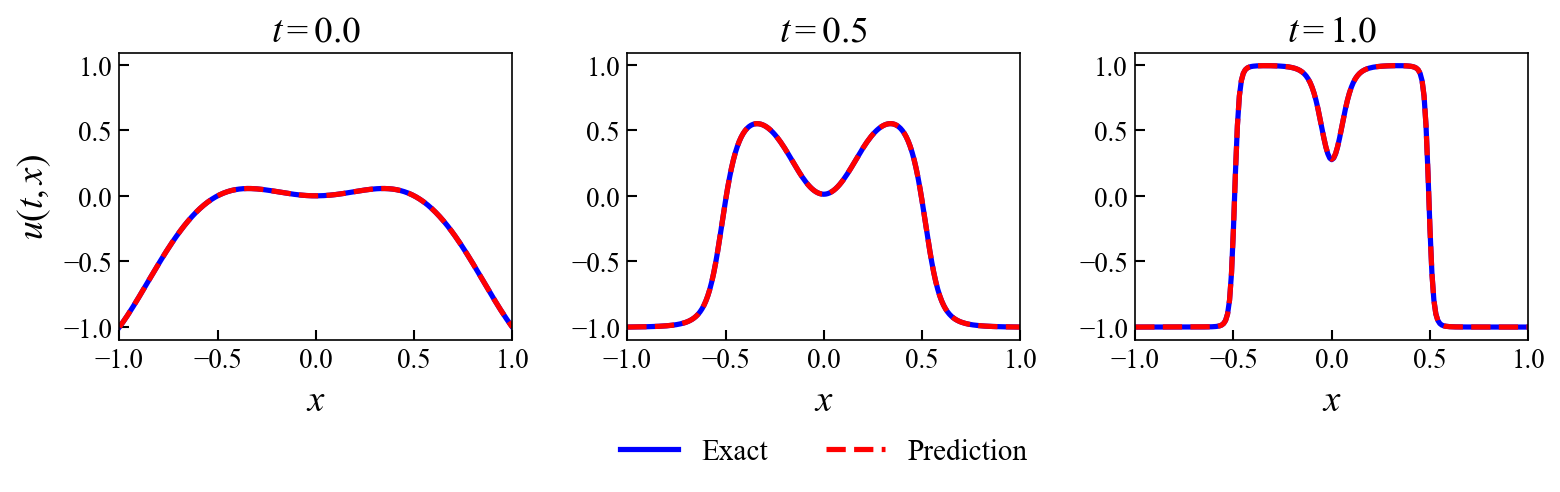}
    \caption{{Allen--Cahn (Dirichlet): cross-sectional
    comparisons of the HC-PINN prediction against the reference solution
    at $t=0$, $t=0.5$, and $t=1.0$.}}
    \label{fig:ac_dirichlet_slices}
\end{figure}

{HC-PINN attains a relative $L^2$ error of
$5.05\times 10^{-4}$ and a relative $L^1$ error of $2.82\times 10^{-4}$,
 slightly lower than  the periodic Case~$I$ error
of $8.29\times 10^{-4}$ reported in Table~\ref{tab:ac}. Figure~\ref{fig:ac_dirichlet}
shows the reference solution, the HC-PINN prediction, and the pointwise
absolute error: the network captures the sharp interface dynamics and
the boundary layers induced by the Dirichlet data without any visible
artifacts at the boundary.}

{This experiment confirms that the stabilizing role of
exact initial-condition enforcement is independent of the periodic
Fourier-feature architecture: the same hard-constraint mechanism,
boundary-compatible choice of $\psi$ and a multiplier $\phi$ that vanishes on $\partial\Omega$, transfers cleanly to the Dirichlet setting. A systematic
study across Neumann, Robin, and mixed boundary types on more complex
geometries is left to future work.}

\section{Conclusion}\label{section:conclude}
We provide a comprehensive study of hard constraints and self-adaptive loss weighting for stabilizing the training of PINNs in solving stiff time-dependent PDEs. Our study includes an ablation analysis across seven representative PDEs, a comparison with state-of-the-art methods, a theoretical investigation, an application to the $2D$ Allen--Cahn equation, and a sensitivity analysis with respect to $m$ (the number of Fourier modes).

The ablation results demonstrate that the hard-constrained transformation plays a dominant role in alleviating the training challenges associated with stiff time-dependent PDEs. In particular, it consistently improves performance across a diverse set of problems, including scalar-valued and vector-valued (including complex-valued) systems. This transformation enhances both training stability and solution accuracy, especially in regimes where conventional PINNs often fail. From a theoretical perspective, our NTK analysis provides insight into this improvement by showing that the hard-constrained formulation mitigates spectral bias and simplifies the training dynamics.

The comparison study further confirms that HC-PINN outperforms existing state-of-the-art approaches on standard benchmarks, highlighting its effectiveness and robustness. Overall, our findings suggest that enforcing initial conditions through hard constraints is not merely beneficial, but essential for reliable PINN training in stiff regimes.

{The reported experiments are predominantly periodic,
with one non-periodic Dirichlet example (Section~\ref{sec:dirichlet})
indicating that the hard-IC mechanism is the dominant stabilizing factor
irrespective of boundary type. A systematic study across Neumann, Robin,
and mixed boundaries on complex geometries is left to future work.} Future work includes extending the hard-constrained framework to more complex two- and three-dimensional PDEs, where additional challenges arise in both architecture design and computational cost. Another important direction is the systematic design of optimal transformations $\psi$ and $\phi$, which govern the embedding of initial and boundary conditions into the network structure, with the goal of further improving performance and generalizability.

\section*{Author contributions}
\noindent \textbf{Baoli Hao:} Conceptualization, Methodology, Software, Investigation, Formal analysis, Visualization, Writing - original draft;
\textbf{Ming Zhong:} Conceptualization, Methodology, Supervision, Writing - review \& editing;
\textbf{Ulisses Braga-Neto:} Conceptualization, Methodology;
\textbf{Chun Liu:} Supervision, Writing - review \& editing;
\textbf{Lifan Wang:} Conceptualization, Writing - review \& editing.
All authors have read and agreed to the published version of the manuscript.

\section*{Acknowledgments}
The authors BH and MZ would like to thank Prof. Yannis Kevrekidis (JHU) for the discussion on the well-conditioning of PDE solution maps, and Prof. George Karniadakis (Brown) for his suggestions on training enhancement methods for PINNs. MZ is supported by NSF-AoF grant \#2225507 and Ralph E. Powe Junior Faculty Enhancement Award from ORAU for FY24.
\section*{Use of Generative-AI tools declaration}
The authors used AI-based language tools to improve clarity and grammar. No AI tools were used for scientific content generation.

%%%%%%%%%%%%%%%%%%%%%%%%%%%%%%%%%%%%%%%%%%%%%%%%%%%%%%
%          7. REFERENCES SECTION
%%%%%%%%%%%%%%%%%%%%%%%%%%%%%%%%%%%%%%%%%%%%%%%%%%%%%%

%       READ THIS SECTION CAREFULLY

% Each of the references below MUST be cited in your article above. Do not include references that are not cited in your article.

% Follow the examples below carefully. We strongly suggest that you copy and paste your reference information directly into our examples.

% List all references in alphabetical order according to the first author's last name.

% Verify each URL works correctly and can be accessed properly. Your URL links should be to reputable websites. The command line for a website link begins with: \url{ }

% Do not add MR or DOI numbers to your references. AIMS production staff will add this information.

% Using BibTex is not recommended but can be handled.

\appendix

\section{Proofs of the Theorems}
\label{section:appendix-proof}

We present the proofs of the two main theoretical results in Section~\ref{sec:theory_full}.
The elementary derivation associated with Theorem~\ref{thm:advection}
is provided directly in Section~\ref{sec:la}.

\section*{Proof of Theorem~\ref{thm:wellposed}
(Consistency and Non-Degeneracy of the Transformation)}

\begin{proof}
By the definition of the hard-constraint transformation,
\[
    u(t,\bx)
    =
    \psi(t,\bx)
    +
    \phi(t,\bx)u_{nn}(t,\bx).
\]
Using
\[
    \psi(0,\bx)=u_0(\bx),
    \qquad
    \phi(0,\bx)=0,
    \qquad
    \bx\in\bar{\Omega},
\]
we immediately obtain
\[
    u(0,\bx)
    =
    \psi(0,\bx)
    +
    \phi(0,\bx)u_{nn}(0,\bx)
    =
    u_0(\bx).
\]
Thus, the prescribed initial condition is satisfied exactly.

Next, differentiating the transformation with respect to time gives
\[
    \partial_tu
    =
    \psi_t
    +
    \phi\partial_tu_{nn}
    +
    \phi_tu_{nn}.
\]
Therefore,
\[
\begin{aligned}
    \partial_tu+\mP[u]-f
    &=
    \psi_t
    +
    \phi\partial_tu_{nn}
    +
    \phi_tu_{nn}
    +
    \mP[\psi+\phi u_{nn}]
    -
    f.
\end{aligned}
\]
Hence,
\[
    \partial_tu+\mP[u]=f
\]
if and only if
\[
    \phi\partial_tu_{nn}
    +
    \mP[\psi+\phi u_{nn}]
    +
    \phi_tu_{nn}
    =
    f-\psi_t.
\]
This proves the equivalence between the original PDE and the transformed PDE.

Finally, evaluating the time derivative of the transformation at $t=0$ and
using $\phi(0,\bx)=0$, we obtain
\[
    u_t(0,\bx)
    =
    \psi_t(0,\bx)
    +
    \phi_t(0,\bx)u_{nn}(0,\bx).
\]
Since
\[
    \phi_t(0,\bx)\neq 0,
    \qquad
    \bx\in\bar{\Omega},
\]
the initial value of the transformed variable is uniquely determined by
\[
    u_{nn}(0,\bx)
    =
    \frac{
        u_t(0,\bx)-\psi_t(0,\bx)
    }{
        \phi_t(0,\bx)
    }.
\]
Therefore, the transformation is non-degenerate at the initial time in the
sense stated in the theorem.
\end{proof}

\subsection*{Proof of Theorem~\ref{thm:hessian} (Hessian Characterization and Bounds)}

\begin{proof}
Since
\[
    u(\by;\btheta)=\h^\top(\by)\btheta,
\]
the loss of the soft-constrained PINN can be written as
\[
\begin{aligned}
    \mathrm{Loss}_s(u(\btheta))
    &=
    \frac{1}{2N}
    \sum_{i=1}^{N}
    \norm{
        \bigl(
            \mL_s[u(\cdot;\btheta)]-f
        \bigr)(\bz_i)
    }^2  \\
    &\qquad
    +
    \frac{1}{2M}
    \sum_{j=1}^{M}
    \norm{
        u(\bzeta_j;\btheta)
        -
        u_0(\bx_j^{IC})
    }^2.
\end{aligned}
\]
Define
\[
    \bA
    =
    \begin{bmatrix}
        \bigl(\mL_s[\h^\top]\bigr)(\bz_1)\\
        \vdots\\
        \bigl(\mL_s[\h^\top]\bigr)(\bz_N)
    \end{bmatrix},
    \qquad
    \bC
    =
    \begin{bmatrix}
        \h^\top(\bzeta_1)\\
        \vdots\\
        \h^\top(\bzeta_M)
    \end{bmatrix},
\]
and
\[
    \f
    =
    \begin{bmatrix}
        f(\bz_1)\\
        \vdots\\
        f(\bz_N)
    \end{bmatrix},
    \qquad
    \bu_0
    =
    \begin{bmatrix}
        u_0(\bx_1^{IC})\\
        \vdots\\
        u_0(\bx_M^{IC})
    \end{bmatrix}.
\]
Then
\[
    \mathrm{Loss}_s(u(\btheta))
    =
    \frac{1}{2N}
    \norm{
        \bA\btheta-\f
    }_{\ell^2(\R^N)}^2
    +
    \frac{1}{2M}
    \norm{
        \bC\btheta-\bu_0
    }_{\ell^2(\R^M)}^2.
\]

For the hard-constrained formulation, recall that
\[
    \widetilde{u}
    =
    \psi+\phi u.
\]
Since
\[
    \mL_h[u]
    =
    \phi\partial_tu
    +
    \mP[\psi+\phi u]
    +
    \phi_tu,
\]
we have
\[
    \mL_h[u]-f+\psi_t
    =
    \mL_s[\psi+\phi u]-f.
\]
Therefore, the hard-constrained loss can equivalently be written as
\[
    \mathrm{Loss}_h(u(\btheta))
    =
    \frac{1}{2N}
    \sum_{i=1}^{N}
    \norm{
        \bigl(
            \mL_s[\psi(\cdot)+\phi(\cdot)u(\cdot;\btheta)]
            -
            f
        \bigr)(\bz_i)
    }^2.
\]

Since $\mP$ is linear and $\phi(t,\bx)=\phi(t)$, the operator
$\mP$ acts only on the spatial variables and hence
\[
    \mP[\phi u]=\phi\mP[u].
\]
Consequently,
\[
\begin{aligned}
    \mL_s[\psi+\phi u]
    &=
    \partial_t(\psi+\phi u)
    +
    \mP[\psi+\phi u] \\[0.2cm]
    &=
    \psi_t
    +
    \phi_tu
    +
    \phi u_t
    +
    \mP[\psi]
    +
    \phi\mP[u] \\[0.2cm]
    &=
    \phi\bigl(u_t+\mP[u]\bigr)
    +
    \phi_tu
    +
    \psi_t
    +
    \mP[\psi] \\[0.2cm]
    &=
    \phi\mL_s[u]
    +
    \phi_tu
    +
    \mL_s[\psi].
\end{aligned}
\]
Define
\[
    \bB
    =
    \begin{bmatrix}
        \h^\top(\bz_1)\\
        \vdots\\
        \h^\top(\bz_N)
    \end{bmatrix},
    \qquad
    \tilde{\f}
    =
    \begin{bmatrix}
        f(\bz_1)-\mL_s[\psi](\bz_1)\\
        \vdots\\
        f(\bz_N)-\mL_s[\psi](\bz_N)
    \end{bmatrix},
\]
and
\[
    \bLambda_{\phi}
    =
    \operatorname{diag}
    \bigl(
        \phi(\bz_1),\ldots,\phi(\bz_N)
    \bigr),
    \qquad
    \bLambda_{\phi_t}
    =
    \operatorname{diag}
    \bigl(
        \phi_t(\bz_1),\ldots,\phi_t(\bz_N)
    \bigr).
\]
It follows that
\[
    \mathrm{Loss}_h(u(\btheta))
    =
    \frac{1}{2N}
    \norm{
        \left(
            \bLambda_{\phi}\bA
            +
            \bLambda_{\phi_t}\bB
        \right)\btheta
        -
        \tilde{\f}
    }_{\ell^2(\R^N)}^2.
\]

The corresponding gradients are
\[
\begin{aligned}
    \nabla_{\btheta}\mathrm{Loss}_s(u(\btheta))
    &=
    \frac{1}{N}
    \left(
        \bA^\top\bA\btheta
        -
        \bA^\top\f
    \right)
    +
    \frac{1}{M}
    \left(
        \bC^\top\bC\btheta
        -
        \bC^\top\bu_0
    \right), \\
    \nabla_{\btheta}\mathrm{Loss}_h(u(\btheta))
    &=
    \frac{1}{N}
    \Bigl[
        \left(
            \bLambda_{\phi}\bA
            +
            \bLambda_{\phi_t}\bB
        \right)^\top
        \left(
            \bLambda_{\phi}\bA
            +
            \bLambda_{\phi_t}\bB
        \right)\btheta \\
    &\hspace{3.8cm}
        -
        \left(
            \bLambda_{\phi}\bA
            +
            \bLambda_{\phi_t}\bB
        \right)^\top
        \tilde{\f}
    \Bigr].
\end{aligned}
\]
Differentiating once more with respect to $\btheta$, we obtain
\[
    H_s(\btheta)
    =
    \frac{1}{N}\bA^\top\bA
    +
    \frac{1}{M}\bC^\top\bC,
\]
and
\[
\begin{aligned}
    H_h(\btheta)
    &=
    \frac{1}{N}
    \left(
        \bLambda_{\phi}\bA
        +
        \bLambda_{\phi_t}\bB
    \right)^\top
    \left(
        \bLambda_{\phi}\bA
        +
        \bLambda_{\phi_t}\bB
    \right) \\
    &=
    \frac{1}{N}
    \Bigl(
        \bA^\top\bLambda_{\phi}^2\bA
        +
        \bA^\top\bLambda_{\phi}\bLambda_{\phi_t}\bB
        +
        \bB^\top\bLambda_{\phi_t}\bLambda_{\phi}\bA
        +
        \bB^\top\bLambda_{\phi_t}^2\bB
    \Bigr).
\end{aligned}
\]

We next derive the Hessian bounds. Here,
$\norm{\mL_s}_2$ denotes an induced discrete operator bound on the sampled
linearized trial space such that
\[
    \norm{\bA}_2
    \leq
    \norm{\mL_s}_2\norm{\bB}_2.
\]
Therefore,
\[
    \norm{\bA^\top\bA}_2
    =
    \norm{\bA}_2^2
    \leq
    \norm{\mL_s}_2^2\norm{\bB}_2^2.
\]
It follows immediately that
\[
\begin{aligned}
    \norm{H_s}_2
    &\leq
    \frac{1}{N}
    \norm{\bA^\top\bA}_2
    +
    \frac{1}{M}
    \norm{\bC^\top\bC}_2 \\
    &\leq
    \frac{\norm{\mL_s}_2^2}{N}
    \norm{\bB}_2^2
    +
    \frac{1}{M}
    \norm{\bC}_2^2.
\end{aligned}
\]

Define
\[
    C_1
    =
    \max_{\bz}\phi^2(\bz),
    \qquad
    C_2
    =
    \max_{\bz}
    \left|
        \phi(\bz)\phi_t(\bz)
    \right|,
    \qquad
    C_3
    =
    \max_{\bz}\phi_t^2(\bz).
\]
Then
\[
\begin{aligned}
    \norm{
        \bA^\top\bLambda_{\phi}^2\bA
    }_2
    &\leq
    C_1\norm{\bA}_2^2
    \leq
    C_1\norm{\mL_s}_2^2\norm{\bB}_2^2, \\[0.2cm]
    \norm{
        \bA^\top
        \bLambda_{\phi}
        \bLambda_{\phi_t}
        \bB
    }_2
    &\leq
    C_2\norm{\bA}_2\norm{\bB}_2
    \leq
    C_2\norm{\mL_s}_2\norm{\bB}_2^2, \\
    \norm{
        \bB^\top
        \bLambda_{\phi_t}
        \bLambda_{\phi}
        \bA
    }_2
    &\leq
    C_2\norm{\bB}_2\norm{\bA}_2
    \leq
    C_2\norm{\mL_s}_2\norm{\bB}_2^2, \\
    \norm{
        \bB^\top\bLambda_{\phi_t}^2\bB
    }_2
    &\leq
    C_3\norm{\bB}_2^2.
\end{aligned}
\]
Combining these estimates yields
\[
\begin{aligned}
    \norm{H_h}_2
    &\leq
    \frac{1}{N}
    \left(
        C_1\norm{\mL_s}_2^2
        +
        2C_2\norm{\mL_s}_2
        +
        C_3
    \right)
    \norm{\bB}_2^2.
\end{aligned}
\]
This completes the proof.
\end{proof}
%\bibliographystyle{plain}
% \bibliography{cas-refs}
% \bibliographystyle{johd}
%\bibliographystyle{original_doc/johd}

\end{document}